\documentclass[12pt]{amsart}
\usepackage{amssymb,amsmath}
\usepackage[UglyObsolete]{diagrams}
\usepackage[hypertex]{hyperref}
\swapnumbers
\newtheorem{theorem}{Theorem}[section]
\newtheorem{corollary}[theorem]{Corollary}

\newtheorem{lemma}[theorem]{Lemma}
\newtheorem{proposition}[theorem]{Proposition}

\theoremstyle{definition}
\newtheorem{example}[theorem]{Example}
\newtheorem*{warning}{Warning}
\newtheorem{convention}[theorem]{Convention}
\newtheorem{remark}[theorem]{Remark}
\newtheorem{definition}[theorem]{Definition}
\newtheorem{exercise}[theorem]{Exercise}
\newtheorem*{formality-theorem}{Formality}

\newcommand{\mb}{\mathbb}
\newcommand{\mc}{\mathcal}
\newcommand{\mf}{\mathfrak}
\newcommand{\ms}{\mathsf}

\newcommand{\abs}[1]{{\lvert #1 \rvert}}

\newcommand{\Der}{\mathop{\mathit{Der}}\nolimits}
\newcommand{\IDer}{\mathop{\mathit{IDer}}\nolimits}
\newcommand{\CoDer}{\mathop{\mathit{CoDer}}\nolimits}
\newcommand{\Diff}{\mathop{\mathit{Diff}}\nolimits}
\newcommand{\CoDiff}{\mathop{\mathit{CoDiff}}\nolimits}
\newcommand{\Span}{\mathop{\mathrm{Span}}\nolimits}

\def\To^#1{\stackrel{#1}{\to}}
\newcommand{\Tto}{\mathop{\longrightarrow}\limits} 
\newcommand{\id}{\mathop{\mathit{id}}\nolimits} 
\newcommand{\Ker}{\mathop{\mathrm{Ker}}\nolimits}
\newcommand{\GL}{\mathop{\mathit{GL}}\nolimits}

\newcommand{\sgn}{\mathop{\mathrm{sgn}}\nolimits}

\newcommand{\Alg}{\mathop{\mathsf{Alg}}\nolimits}

\newcommand{\Ass}{\mathop{\mathsf{Ass}}\nolimits}

\newcommand{\Com}{\mathop{\mathsf{Com}}\nolimits}
\newcommand{\Lie}{\mathop{\mathsf{Lie}}\nolimits}
\newcommand{\Linf}{\mathop{\mathsf{L}_\infty}\nolimits}

\def\bfk{\mathbf {k}} \def\Hoch{{\rm Hoch}} 
\def\Lin{{\it Lin}}\def\otexp#1#2{{#1^{\otimes #2}}}
\def\ot{\otimes} \def\Coass{{\it Coass\/}}
\def\Asss{{\it Ass\/}}
\def\susp{\mathop{\uparrow}}
\def\desusp{\mathop{\downarrow}}
\def\dessp#1{\hbox{$\downarrow \hskip -.2em  #1$}}
\def\sssp#1{\hbox{$\uparrow \hskip -.2em  #1$}}
\def\ext{\mbox{\large$\land$}} \def\bvee{\mbox{\large$\vee$}}
\def\squeezedldots{{\mbox{$. \hskip -1pt . \hskip -1pt .$}}}

\def\cT{{{}^cT}} \def\cext{{{}^c\ext}} \def\Gg{{\mbox {\rm G}(\mf{g})}}
\def\MC{{\rm MC}} \def\Def{\mf{Def}}
\def\catLinfty{{\sf L}_\infty}
\def\minibullet{\mbox {\scriptsize $\bullet$}}

\begin{document}

\bibliographystyle{plain}
\baselineskip18pt plus 1pt minus 1pt
\parskip3pt plus 1pt minus .5pt

\title[Deformation Theory]{Deformation Theory
\\
\rule{0pt}{1.3em}
(lecture notes)}
\author[Doubek, Markl, Zima]{M.~Doubek, M.~Markl and P.~Zima}
\thanks{The second author was supported 
   by the grant GA \v CR 201/05/2117 and by
   the Academy of Sciences of the Czech Republic,
   Institutional Research Plan No.~AV0Z10190503}

\address{\hbox{\bf Corresponding author:} M.~Markl: Math. 
         Inst. of the Academy, {\v Z}itn{\'a} 25,
         115 67 Prague~1, The Czech Republic}
\email{martindoubek@seznam.cz, markl@math.cas.cz}

\maketitle

\begin{center}
{\rm
Notes, taken by Martin~Doubek and Petr Zima, 
from a course given by Martin~Markl at~the~Charles University,
Prague, in the Summer semester 2006.}
\end{center}

\begin{abstract}
First three sections of this overview paper
cover classical topics of deformation theory of associative algebras 
and necessary background material. We then 
analyze algebraic structures of the Hochschild cohomology and describe
the relation between deformations and solutions of the corresponding
Maurer-Cartan equation. In Section~\ref{s7}  we generalize the
Maurer-Cartan equation to strongly homotopy Lie
algebras and prove the homotopy invariance of the moduli space of
solutions of this equation. In the last section we indicate the main
ideas of Kontsevich's proof of the existence of deformation
quantization of Poisson manifolds. 
\end{abstract}

\vskip 3mm
\noindent 
{\bf Table of content:} \ref{s1}.  
                     Algebras and modules
                     -- p.~\pageref{s1}
\hfill\break\noindent 
\hphantom{{\bf Table of content:\hskip .5mm}}  \ref{s2}.
                     Cohomology -- p.~\pageref{s2}
\hfill\break\noindent 
\hphantom{{\bf Table of content:\hskip .5mm}}  \ref{s3}.
                     Classical deformation theory -- p.~\pageref{s3}
\hfill\break\noindent 
\hphantom{{\bf Table of content:\hskip .5mm}}  \ref{s4}. 
                    Structures of (co)associative (co)algebras 
                    -- p.~\pageref{s4}
\hfill\break\noindent 
\hphantom{{\bf Table of content:\hskip .5mm}}  \ref{s5}.
                      dg-Lie algebras  and the
                     Maurer-Cartan equation -- p.~\pageref{s5}
\hfill\break\noindent 
\hphantom{{\bf Table of content:\hskip .5mm}}  \ref{s7}.
         $L_\infty$-algebras and the Maurer-Cartan equation  -- p.~\pageref{s7}
\hfill\break\noindent 
\hphantom{{\bf Table of content:\hskip .5mm}}  \ref{s8}.
          Homotopy invariance of the Maurer-Cartan equation  -- p.~\pageref{s8}
\hfill\break\noindent 
\hphantom{{\bf Table of content:\hskip .5mm}}  \ref{s9}.
          Deformation quantization of Poisson manifolds  
          -- p.~\pageref{s9}

\vskip .3em
\noindent 
{\bf Conventions.}
All algebraic objects
will be considered over a fixed field $\mathbf {k}$ of characteristic zero.
The symbol $\otimes$ will denote the tensor product over $\bfk$. We
will sometimes use the same symbol for both an 
algebra and its underlying space.

\vskip .3em
\noindent 
{\bf Acknowledgement.} We would like to thank Dietrich Burde for
useful comments on a~preliminary version of this paper. We are also
indebted to Ezra Getzler for turning our attention to a remarkable
paper~\cite{dauria-fre:NP82}. Also suggestions of M.~Goze and
E.~Remm were very helpful.

\section{Algebras and modules}
\label{s1}

In this section we investigate modules (where module  means rather 
a~bimodule than a~one-sided module) 
over various types of algebras.

\begin{example}\label{e1} -- 
The category $\Ass$ of associative algebras.\\
An associative algebra is a $\mathbf {k}$-vector space $A$
with a bilinear multiplication $A \otimes A \to A$ satisfying
\begin{eqnarray*}
a(bc) = (ab)c,
\qquad \textrm{for all } a,b,c \in A.
\end{eqnarray*}
Observe that at this moment we do not assume the existence of a unit
$1 \in A$.

What we understand by a module over an associative algebra 
is in fact a bimodule, i.e.~a~vector space
$M$ equipped with multiplications (``actions'') 
by elements of $A$ from both sides, subject to the axioms
\begin{eqnarray*}
a(bm) & = & (ab)m, \\
a(mb) & = & (am)b, \\
m(ab) & = & (ma)b,
\qquad \textrm{for all } m \in M,\, a,b \in A.
\end{eqnarray*}
\end{example}

\begin{example} --
The category $\Com$ of commutative associative algebras.\\
In this case left modules, right modules and bimodules coincide.
In addition to the axioms in $\Ass$ we require the commutativity
\[
ab = ba, \qquad \textrm{for all } a,b \in A,
\]
and for a module
\begin{eqnarray*}
ma = am,
\qquad \textrm{for all } m \in M,\, a \in A.
\end{eqnarray*}
\end{example}

\begin{example}\label{e3} --
The category $\Lie$ of Lie algebras.\\
The bilinear bracket $[-,-] : L \otimes L \to L$ of a Lie algebra $L$ 
is anticommutative and satisfies the Jacobi identity, that is
\begin{eqnarray*}
[a,b] & = & -[b,a], \mbox { and}
\\ \relax
[a,[b,c]] + [b,[c,a]] + [c,[a,b]] & = & 0,
\qquad \textrm{for all } a,b,c \in L.
\end{eqnarray*}
A left module (also called a representation) $M$ of $L$ satisfies 
the standard axiom
\[
a(bm) - b(am) = [a,b]m,
\qquad \textrm{for all } m \in M,\, a,b \in L.
\]
Given a left module $M$ as above, 
one can canonically turn  it into a right module by setting
$ma := -am$. Denoting these actions of $L$ by the bracket, one can rewrite
the axioms as
\begin{eqnarray*}
[a,m] & = & -[m,a], \mbox { and}
\\ \relax
[a,[b,m]] + [b,[m,a]] + [m,[a,b]] & = & 0,
\qquad \textrm{for all } m \in M,\, a,b \in L.
\end{eqnarray*}
\end{example}

Examples~\ref{e1}--\ref{e3} indicate how axioms of algebras induce, by
replacing one instance of an algebra variable by a module variable,
axioms for the corresponding modules.  In the rest of this section we
formalize, following~\cite{quillen}, this recipe. The standard
definitions below can be found for example in~\cite{maclane:71}.

\begin{definition}
The \emph{product} in a 
category $\mathsf{C}$ is the limit of a discrete diagram.
The \emph{terminal object} of $\mathsf{C}$ is the limit of an empty
diagram, or equivalently, an object $T$ such that 
for every $X \in \mathsf{C}$ there
exists a unique morphism $X \to T$.
\end{definition}

\begin{remark}
The product of any object $X$ with the terminal object 
$T$ is naturally isomorphic to $X$,
\[
X \times T \cong X \cong T \times X.
\]
\end{remark}

\begin{remark}
\label{R1}
It follows from the universal property of the product that there
exists the \emph{swapping morphism} $X\times X\To^s X\times X$ making
the diagram
\begin{diagram}
X\times X & \rTo^{p_1} & X \\
\dTo<{p_2} & \rdDashto^{s} & \uTo>{p_2} \\
X & \lTo_{p_1} & \hphantom{.}X\times X,
\end{diagram}
in which $p_1$ (resp.~$p_2$) is the projection onto the first
(resp.~second) factor, commutative.
\end{remark}

\begin{example}
In the category of $A$-bimodules, the product $M_1 \times M_2$ is the ordinary
direct sum $M_1 \oplus M_2$.
The terminal object is the trivial module~$0$.
\end{example}

\begin{definition}
A category $\mathsf{C}$ has \emph{finite products}, if every \emph{finite}
discrete diagram has a limit in $\mathsf{C}$.
\end{definition}

By~\cite[Proposition~5.1]{maclane:71}, $\mathsf{C}$ has finite limits
if and only if it has a terminal object and products of pairs of objects.

\begin{definition}
Let $\mathsf{C}$ be a category, $A\in\mathsf{C}$. The 
\emph{comma category} (also called the \emph{slice category}) 
$\mathsf{C}/A$ is the category whose

-- objects $(X,\pi)$ are  $\mathsf{C}$-morphisms $X\To^\pi A$, 
$X\in\mathsf{C}$, and 

-- morphisms $(X',\pi') \To^{f} (X'',\pi'')$ are commutative
diagrams of $\mathsf{C}$-morphisms:
\begin{diagram}
X' & \rTo^f & X'' \\
\dTo<{\pi'} & & \dTo>{\pi''} & \\
A & \rEq_{\id_A} & \hphantom{.}A.
\end{diagram}
\end{definition}

\begin{definition}
The \emph{fibered product} (or \emph{pullback}) of morphisms
$X_1\To^{f_1} A$ and $X_2\To^{f_2} A$ in $\mathsf{C}$ is the limit
$D$ (together with morphisms $D\To^{p_1} X_1,\ D\To^{p_2} X_2$) of the lower
right corner of the digram:
\begin{diagram}
D & \rDashto^{p_1} & X_1 \\
\dDashto<{p_2} & & \dTo>{f_1} \\
X_2 & \rTo_{f_2} & \hphantom{.}A.
\end{diagram}
In the above situation one sometimes writes $D = X_1 \times_A X_2$.
\end{definition}

\begin{proposition}
If $\mathsf{C}$ has fibered products then $\mathsf{C}/A$ has finite
products.
\end{proposition}

\begin{proof} A straightforward verification.
The identity morphism $(A,\id_A)$ is clearly the terminal object
of $\mathsf{C}/A$.

Let $(X_1,\pi_1)$ and $(X_2,\pi_2)$ be objects of $\mathsf{C}/A$. 
By assumption, there exists the fibered product
\begin{equation}
\label{EE1}
\begin{diagram}
D & \rDashto^{p_1} & X_1 \\
\dDashto<{p_2} & \rdDashto^\delta & \dTo> {\pi_1} \\
X_2 & \rTo_{\pi_2} & A
\end{diagram}
\end{equation}
in $\mathsf{C}$. In the above diagram, of course, $\delta := \pi_1 p_1
= \pi_2 p_2$.
The maps $p_1 : D \to X_1$ and $p_2 : D \to X_2$ of the above diagram
define morphisms (denoted by the same symbols) $p_1: (D,\delta)
\to (X_1,\pi_1)$ and $p_2: (D,\delta) \to (X_2,\pi_2)$ in $\mathsf{C}/A$. 
The universal property of the
pullback~(\ref{EE1}) implies that the object $(D,\delta)$ with the
projections $(p_1,p_2)$ is the product of $(X_1,\pi_1) \times
(X_2,\pi_2)$ in~$\mathsf{C}/A$.       
\end{proof}

One may express the conclusion of the above proof by
\begin{equation}
\label{E187}
(X_1,\pi_1) \times (X_2,\pi_2) = X_1 \times_A X_2,
\end{equation}  
but one must be aware that the left side lives in $\mathsf{C}/A$
while the right one in $\mathsf{C}$, therefore~(\ref{E187}) has only
a symbolical meaning. 

\begin{example}
In $\Ass$, the fibered product of morphisms $B_1\To^{f_1} A,\ B_2\To^{f_2}A$
is the subalgebra
\begin{equation}
\label{E23}
B_1\times_A B_2 = \{(b_1,b_2)\ |\ f_1(b_1) =
f_2(b_2)\} \subseteq B_1 \oplus B_2
\end{equation}
together with the restricted projections. Hence for any algebra $A \in \Ass$,
the comma category $\Ass/A$ has finite products.
\end{example}

\begin{definition}
Let $\mathsf{C}$ be a category with 
finite products and $T$ its terminal object.
An~\emph{abelian group object} in $\mathsf{C}$ is a 
quadruple $(G,G\times G\To^{\mu} G,G \To^\eta G, T\To^e G)$
of objects and morphisms of $\mathsf{C}$ such that following
diagrams commute:

-- the \emph{associativity} $\mu$:
\begin{diagram}
G\times G\times G & \rTo^{\mu\times\id_G} & G\times G \\
\dTo<{\id_G\times\mu} & & \dTo>\mu \\
G\times G & \rTo_\mu & \hphantom{,}G,
\end{diagram}

-- the \emph{commutativity} of $\mu$ (with $s$ the swapping 
morphism of Remark~\ref{R1}):
\begin{diagram}
G\times G & & \rTo^s & & G\times G \\
 & \rdTo<\mu & & \ldTo>\mu & \\
  & & G & &
\end{diagram}

-- the \emph{neutrality} of $e$:
\begin{diagram}
T\times G & \rTo^{e\times\id_G} & G\times G &
\lTo^{\id_G\times e} & G\times T \\
 & \rdEq>\cong & \dTo>\mu & \ldEq<\cong & \\
 & & G & &
\end{diagram}

-- the diagram saying that $\eta$ is a two-sided \emph{inverse} for the
   multiplication $\mu$:
\begin{diagram}
G & \rTo^{\eta \times {\id}_G} & G \times G \\
\dTo<{{\id}_G \times \eta} & \rdTo& \dTo>{\mu} \\
G \times G & \rTo_{\mu} & \hphantom{.}G,
\end{diagram}
in which the diagonal map is the composition $G \to T \stackrel e\to G$.

Maps $\mu$, $\eta$ and $e$ above are called the \emph{multiplication},
the \emph{inverse} and the \emph{unit} of the abelian group
structure, respectively.

\emph{Morphisms} of abelian group objects $(G',\mu',\eta',e')
\To^f(G'',\mu'',\eta'',e'')$ are
morphisms $G' \To^f G''$ in $\mathsf{C}$ which preserve all structure
operations. In terms of diagrams this means that
\begin{diagram}
G'\times G' & \rTo^{f\times f} & G''\times G'' & \qquad &
G' & \rTo^{f} & G''& \qquad &
T & \rEq^{\id_T} & T 
\\
\dTo<{\mu'} & & \dTo>{\mu''} & &
\dTo<{\eta'} & & \dTo>{\eta''} & &
\dTo<{e'} & & \dTo>{e''} 
\\
G & \rTo_f & G' & &
G & \rTo_f & G' & &
G & \rTo_f & G'
\end{diagram}
commute. The category of abelian group objects of $\mathsf{C}$ will be denoted
$\mathsf{C}_{ab}$.
\end{definition}

Let $\Alg$ be any of the examples of categories of algebras considered
above and $A \in \Alg$. It turns out that the
category $(\Alg/A)_{ab}$ is precisely the corresponding category
of $A$-modules.
To verify this for associative algebras, we identify, in
Proposition~\ref{p1} below, objects of
$(\Ass/A)_{ab}$ with trivial extensions in the sense of:

\begin{definition}
Let $A$ be an associative algebra and $M$ an $A$-module. The \emph{trivial
extension} of $A$ by $M$ is the associative algebra $A\oplus M$ with
the multiplication given by 
\[
(a,m)(b,n) = (ab, an + mb),\
a,b \in A \mbox { and } m,n \in M.
\]
\end{definition}

\begin{proposition}
\label{p1}
The category $(\Ass/A)_{ab}$ is isomorphic to the category of trivial
extensions of $A$.
\end{proposition}

\begin{proof}
Let $M$ be an $A$-module and $A \oplus M$ the corresponding trivial 
extension. Then $A\oplus M$ with the projection 
$A\oplus M\To^{\pi_A} A$ determines  
an object $G$  of $\Ass/A$ and, by~(\ref{E187}) and~(\ref{E23}), $G \times
G = (A \oplus M \oplus M \To^{\pi_A} A)$. It is clear that
$\mu : G \times G \to G$ given by $\mu(a,m_1,m_2) := (a,m_1+m_2)$, $e$
the inclusion $A \hookrightarrow A \oplus M$ and $\eta : G \to G$
defined by $\eta(a,m) := (a,-m)$ make $G$ an abelian group object in
$(\Ass/A)_{ab}$. 

On the other hand, let $((B,\pi),\mu,\eta,e)$ be an abelian group object in
$\Ass/A$. The diagram
\begin{diagram}
A & \rTo^e & B \\
 & \rdEq<{\id_A} & \dTo^\pi \\
 & & A
\end{diagram}
for the neutral element 
says that $\pi$ is a retraction. Therefore one may
identify the algebra $A$ with its image $e(A)$, which is a subalgebra
of $B$. Define $M: =\Ker\pi$ so that there is a vector spaces isomorphism
$B=A\oplus M$ determined by the inclusion $e : A
\hookrightarrow B$ and its retraction~$\pi$. Since $M$ is an ideal in $B$, the
algebra $A$ acts on $M$ from both sides. Obviously, 
$M$ with these actions is an $A$-bimodule, the
bimodule axioms following from the associativity of $B$ as in Example~\ref{e1}.
It remains to show that $m'm''=0$ for all $m',m''\in M$ which would imply 
that $B$ is a trivial extension of $A$. Let us introduce the following
notation. 

For a morphism $f: (B',\pi') \to (B'',\pi'')$ of $\bfk$-splitting objects of 
$\Ass/A$ (i.e.~objects with specific $\bfk$-vector space isomorphisms
$B' \cong  A
\oplus M'$ and  $B'' \cong A \oplus M''$ such that $\pi'$ and $\pi''$
are the projections on the first summand) 
we denote by $\tilde{f} : M' \to M''$ the restriction $f|_{M'}$ followed
by the projection $B'' \To^{\pi'} M''$. We call $\tilde{f}$ the {\em
reduction\/} of $f$. Clearly, for every diagram of splitting objects in
$\Ass/A$ there is the corresponding diagram of reductions in $\Ass$.

The fibered product $(A\oplus M,\pi)\times
(A\oplus M,\pi)$ in $\Ass/A$ is isomorphic to $A \oplus M \oplus
M$ with the multiplication
\[
(a',m_1',m_2')(a'',m_1'',m_2'') = (a'a'',a'm_1'' + m_1'a''+ m_1'm_1'',
a'm_2'' + m_2'a''+ m_2'm_2'').
\]
The neutrality of $e$ implies the following diagram of reductions
\begin{diagram}
0\oplus M & \rTo^{\tilde{e}\times\id_M} & M\oplus M &
\lTo^{\id_M\times\tilde{e}} & M\oplus 0 \\
 & \rdEq^\cong & \dTo>{\tilde{\mu}} & \ldEq^\cong & \\
 & & M & &
\end{diagram}
which in turn implies
\[
\tilde{\mu}(0,m) = \tilde{\mu}(m,0) = m,
\qquad \textrm{for all } m \in M.
\]
Since $\mu$ is a morphism in $\Ass$, it preserves the 
multiplication and so does
its reduction $\tilde{\mu}$. We finally obtain
\[
m'\cdot m''=\tilde\mu(m',0)\cdot \tilde\mu(0,m'')=\tilde\mu((m',0)\cdot(0,m''))
=\tilde\mu(m'\cdot 0,0\cdot m'')=0.
\]
This finishes the proof.
\end{proof}

We have shown that objects of $(\Ass/A)_{ab}$ are precisely trivial
extensions of $A$. Since there is an obvious equivalence between modules
and trivial extensions, we obtain:

\begin{theorem}
\label{t1}
The category $(\Ass/A)_{ab}$ is isomorphic to the category of $A$-modules.
\end{theorem}

\begin{exercise}
Prove analogous statements also 
for $(\Com/A)_{ab}$ and $(\Lie/L)_{ab}$.
\end{exercise}

\begin{exercise}
The only property of abelian group objects used in our proof of
Proposition~\ref{p1} was the existence of a neutral element for the
multiplication. In fact, by analyzing our arguments we conclude that
in $\Ass/A$, every object with a multiplication and a neutral element
(i.e.~a \emph{monoid} in $\Ass/A$) is an abelian group
object. Is this statement true in any comma category? If not, what special
property of $\Ass/A$ makes it hold in this particular category?
\end{exercise}

\section{Cohomology}
\label{s2}

Let $A$ be an algebra, $M$ an $A$-module. There are the following
approaches to the ``cohomology of $A$ with coefficients in $M$.''
\begin{enumerate}
\item[(1)] \emph{Abelian cohomology} 
      defined as $H^*(\Lin(R_*,M))$, where $R_*$ is a 
      resolution of $A$ in the category of $A$-modules.
\item[(2)] \emph{Non-abelian cohomology} defined as 
      $H^*(\Der(\mc{F}_*,M))$, where $\mc{F}_*$ is a resolution of $A$ in
      the category of algebras and $\Der(-,M)$ denotes the space of
      derivations with coefficients in $M$.
\item[(3)] \emph{Deformation cohomology} which is the subject of this
      note. 
\end{enumerate}

The adjective {\em (non)-abelian\/} reminds us that~(1) is a derived functor
in the abelian category of modules while cohomology~(2) is a
derived functor in the non-abelian category of algebras. Construction~(1)
belongs entirely into classical homological
algebra~\cite{maclane:homology}, but (2)~requires Quillen's theory of
closed model categories~\cite{quillen:LNM43}. Recall that in this note
we work over a field of characteristics $0$, over the integers one
should take in~(2) a suitable simplicial resolution~\cite{andre:LNM}.
Let us indicate the meaning of deformation cohomology in the
case of associative algebras.

Let $V=\Span\{e_1,\ldots,e_d\}$ be a 
$d$-dimensional $\bfk$-vector space.\label{p8}
Denote $\Asss(V)$ the set of all associative algebra structures on $V$.
Such a structure is determined by a bilinear map $\mu:V\otimes V\to V$.
Relying on Einstein's convention, 
we write 
$\mu(e_i,e_j)=\Gamma_{ij}^le_l$ for some scalars
$\Gamma_{ij}^l\in \bfk$. The associativity 
$\mu(e_i,\mu(e_j,e_k))=\mu(\mu(e_i,e_j),e_k)$ of $\mu$ can then be
expressed~as
\[
\Gamma_{il}^r\Gamma_{jk}^l=\Gamma_{ij}^l\Gamma_{lk}^r,
\qquad i,j,k,r=1,\ldots,d.
\]
These $d^4$ polynomial equations define an affine algebraic variety, 
which is just another way to view $\Asss(V)$, since every point of 
this variety corresponds to an associative algebra structure on $V$.
We call $\Asss(V)$ the \emph{variety of structure constants} of
associative algebras.

The next step is to consider the quotient $\Asss(V)/\GL(V)$ of
$\Asss(V)$ 
modulo the action of the general linear group $\GL(V)$ recalled in
formula~(\ref{E732}) below.
However, $\Asss(V)/\GL(V)$ is no longer an affine variety, but
only a (possibly singular) algebraic stack (in the sense of Grothendieck).
One can remove singularities by replacing $\Asss(V)$ by a smooth
dg-scheme~$\mc{M}$.
Deformation cohomology is then the cohomology 
of the tangent space of this smooth
dg-scheme~\cite{fontanine-kapranov,fontanine-kapranov:JAMS}.

Still more general approach to deformation cohomology is based on
considering a given category of algebras as the category of
algebras over a certain PROP $\sf P$ and defining the deformation cohomology
using a resolution of $\sf P$ in the category of
PROPs~\cite{kontsevich-soibelman,markl:JPAA96,markl:ib}. When $\sf P$ is a
Koszul quadratic operad, we get the \emph{operadic cohomology} whose
relation to deformations was studied in~\cite{balavoine:CM97}. There
is also an approach to deformations based on 
\emph{triples}~\cite{fox:JPAA93}.  

For associative algebras all the above
approaches give the classical
Hochschild cohomology (formula~3.2 of~\cite[\S X.3]{maclane:homology}):

\begin{definition}
\label{D5}
The \emph{Hochschild cohomology} of an 
associative algebra $A$ with coefficients in an $A$-module $M$ 
is the cohomology of the complex:
\[
0 \Tto M 
\Tto^{\raisebox{-.3em}{\rule{0pt}{0pt}}\delta_\Hoch} C^1_\Hoch(A,M) 
\Tto^{\raisebox{-.3em}{\rule{0pt}{0pt}}\delta_\Hoch} \cdots 
\Tto^{\raisebox{-.3em}{\rule{0pt}{0pt}}\delta_\Hoch} 
C^n_\Hoch(A,M) 
\Tto^{\raisebox{-.3em}{\rule{0pt}{0pt}}\delta_\Hoch}\cdots
\]
in which $C^n_\Hoch(A,M):=\Lin(A^{\otimes n},M)$, the space of
$n$-multilinear maps from $A$ to $M$. 
The coboundary $\delta = \delta_\Hoch : C^n_\Hoch(A,M)
\to C^{n+1}_\Hoch(A,M)$ is defined by
\begin{eqnarray*}
\delta_\Hoch f(a_0\otimes\ldots\otimes a_n)&:=&
(-1)^{n+1} a_0f(a_1\otimes\ldots\otimes a_n)
+f(a_0\otimes\ldots\otimes a_{n-1})a_n
\\
&& +\sum_{i=0}^{n-1}(-1)^{i+n}f(a_0\otimes\ldots
\otimes a_ia_{i+1}\otimes\ldots\otimes a_n),
\end{eqnarray*}
for $a_i\in A$.
Denote $H^n_\Hoch(A,M):= H^n(C^*_\Hoch(A,M),\delta)$.
\end{definition}

\begin{exercise} 
Prove that
$\delta^2_\Hoch=0$.
\end{exercise}

\begin{example}
A simple computation shows that

-- $H^0_\Hoch(A,M)=\{m\in M\ |\ am-ma=0 \textrm{ for all }a\in A\}$,

-- $H^1_\Hoch(A,M)=\Der(A,M)/\IDer(A,M)$, where $\IDer(A,M)$
   denotes the subspace of internal derivations, i.e.~derivations
   of the form $\vartheta_m(a) = am - ma$ for $a \in A$ and some fixed
   $m\in M$.
   Slightly more difficult is to prove that  

-- $H^2_\Hoch(A,M)$ is the space of isomorphism classes of 
   singular extensions of $A$ by $M$~\cite[Theorem X.3.1]{maclane:homology}.
\end{example}

\section{Classical deformation theory}
\label{s3}

As everywhere in this note, we work over a field $\bfk$ of characteristics
zero and $\otimes$ denotes the tensor product over $\bfk$. 
By a {\em ring\/} we will mean a commutative associative $\bfk$-algebra. 
Let us start with necessary preliminary notions.

\begin{definition}
Let $R$ be a ring with unit $e$ and $\omega:\bfk \to R$ the
homomorphism given by $\omega(1): =e$.
A homomorphism $\epsilon:R\to \bfk$ 
is an \emph{augmentation} of $R$ if $\epsilon\omega=\id_\bfk$ or,
diagrammatically,
\begin{diagram}
R & \rTo^\epsilon & \bfk \\
\uTo<\omega & \ruTo>\id & \\
\hphantom{.}\bfk. & & \\
\end{diagram}
The subspace
$\overline{R}:=\Ker\epsilon$ is called the \emph{augmentation ideal}
of $R$. The \emph{indecomposables} of the augmented ring $R$ are
defined as the quotient $Q(R):=\overline{R}/\overline{R}^2$.
\end{definition}

\begin{example}
The unital ring $\bfk[[t]]$ of formal power series with
coefficients in $\bfk$ is augmented, with
augmentation $\epsilon: \bfk[[t]]\to \bfk$ given by
$\epsilon(\sum_{i\in\mb{N}_0}a_it^i):=a_0$. The unital ring $\bfk[t]$
of polynomials with coefficients in $\bfk$ is augmented 
by $\epsilon(f) := f(0)$, for $f \in \bfk[t]$. The truncated
polynomial rings $\bfk[t]/(t^n)$, $n \geq 1$, are also augmented, with
the augmentation induced by the augmentation of $\bfk[t]$.
\end{example}

\begin{example}
Recall that the group ring $\bfk[G]$ of a finite group $G$ with unit $e$ is the
space of all formal linear combinations $\sum_{g\in G}a_gg$, $a_g \in
\bfk$, with the multiplication
\[
\textstyle
(\sum_{g\in G}a'_gg)(\sum_{g\in G}a''_gg) := 
\sum_{g\in G}\sum_{uv = g} a'_{u}a''_vg 
\]
and unit $1 e$. The ring $\bfk[G]$  is augmented by 
$\epsilon: 
\bfk[G]\to \bfk$ given as 
\[
\textstyle
\epsilon(\sum_{g\in G}a_gg):=\sum_{g\in G}a_g.
\]
\end{example}

\begin{example}
A rather trivial example of a ring that does not admit an augmentation
is provided by any proper extension $K \supsetneq \bfk$ of $\bfk$. If
an augmentation $\epsilon : K \to \bfk$ exists, then
$\Ker\epsilon$ is, as an ideal in a field, trivial, which implies that
$\epsilon$ is injective, which would imply that $K = \bfk$
contradicting the assumption $K \neq \bfk$.
\end{example}

\begin{exercise}
If $\sqrt {-1} \not\in \bfk$, then $\bfk[x]/(x^2+1)$ admits no augmentation.
\end{exercise}

In the rest of this section, $R$ will be an augmented unital ring with an
augmentation \hbox{$\epsilon: R \to \bfk$} and the unit map $\omega :
\bfk \to R$.
By a module we will understand a {\em left\/} module.

\begin{remark}
\label{r13}
A unital augmented ring $R$ is a $\bfk$-bimodule, with the bimodule
structure induced by the unit map $\omega$ in the obvious
manner. Likewise, $\bfk$ is an $R$ bimodule, with the structure
induced by $\epsilon$.  If $V$ is a $\bfk$-module, then $R\otimes V$
is an $R$-module, with the action $r' (r''\otimes v):=r'r''\otimes v$,
for $r',r'' \in R$ and $v \in V$.
\end{remark}

\begin{definition}
\label{d51}
Let $V$ be a $\bfk$-vector space and $R$ a unital $\bfk$-ring.
The \emph{free $R$-module} generated by $V$ 
is an $R$-module $R\langle V \rangle$ together with a 
$\bfk$-linear map $\iota: V\to R\langle V \rangle$ 
with the property that for every $R$-module $W$ and a $\bfk$-linear 
map $V\To^\varphi W$, 
there exists a unique $R$-linear map $\Phi: R\langle V \rangle \to W$ 
such that the following diagram commutes:
\begin{diagram}
V & \rTo^\iota & R\langle V \rangle \\
& \rdTo_\varphi & \dDashto_\Phi \\
& & \hphantom{.}W. \\
\end{diagram}
\end{definition}

This universal property determines 
the free module $R\langle V \rangle$ 
uniquely up to isomorphism. A~concrete model is
provided by the $R$-module $R\otimes V$ recalled in Remark~\ref{r13}.

\begin{definition}
Let $W$ be an $R$-module. The
\emph{reduction} of $W$ is the $\bfk$-module 
$\overline{W}:=\bfk\otimes_R W$, with the $\bfk$-action given
by $k'(k''\otimes_R w):=k'k''\otimes_R w$, for $k',k''
\in \bfk$ and $w \in W$.
\end{definition}

One clearly has $\bfk$-module isomorphisms
$\overline{W} \cong W/\overline{R}W$ and $\overline{R
\langle V \rangle} \cong V$. The
reduction clearly defines a functor from the category of $R$-modules
to the category of $\bfk$-modules.

\begin{proposition}
\label{p3}
If $B$ is an associative $R$-algebra, 
then the reduction $\overline{B}$  is a  $\bfk$-algebra, with the structure
induced by the algebra structure of $B$.
\end{proposition}

\begin{proof}
Since $\overline{B}\simeq B/\overline{R}B$, it suffices to verify 
that $\overline{R}B$ is a two-sided ideal in $B$. But this is simple. 
For $r\in\overline{R}$, $b',b''\in B$ one sees that 
$\mu(rb',b'')=r\mu(b',b'')\in\overline{R}B$, which shows that $\mu
(\overline RB,B) \subset \overline RB$.
The right multiplication by elements of $\overline{R}B$ is discussed similarly.
\end{proof}

\begin{definition}
\label{d56}
Let $A$ be an associative $\bfk$-algebra and $R$ an augmented unital ring. 
An \emph{$R$-deformation} of $A$ is an associative 
$R$-algebra $B$ together with a
$\bfk$-algebra isomorphism $\alpha : \overline{B}
\to A$.

Two $R$-deformations $(B',\overline{B}' \stackrel{\alpha'}{\to} A)$
and $(B'',\overline{B}'' \stackrel{\alpha''}{\to} A)$ of $A$
are  \emph{equivalent} if there exists an 
$R$-algebra isomorphism $\phi: B'\to B''$ such that $\overline {\phi} =
\alpha''^{-1}\circ \alpha'$.
\end{definition}

There is probably not much to be said about $R$-deformations without
additional assumptions on the $R$-module $B$. In this note we
assume that \emph{$B$ is a free $R$-module} or, equivalently,
that
\begin{equation}
\label{Eq87}
\hskip 2cm
B\cong R\otimes A\ \mbox {(isomorphism of $R$-modules)}.
\end{equation}
The above isomorphism identifies $A$
with the $\bfk$-linear subspace $1\ot A$ of $B$ and $A \ot A$ with the
$\bfk$-linear subspace $(1 \ot A)\ot (1\ot A)$ of $B \ot B$.

Another assumption frequently used in algebraic
geometry~\cite[Section~III.\S9]{hartshorne:book}  
is that the
$R$-module $B$ is flat which, by definition, means that the functor $B
\otimes_R -$ is left exact. One then speaks about \emph{flat deformations}. 

In what follows, $R$ will be either a power series
ring $\bfk[[t]]$ or a truncation of the polynomial ring $\bfk[t]$ by
an ideal generated by a power of $t$. All these rings are local
Noetherian rings therefore a finitely generated $R$-module is flat if
and only if it is free (see Exercise~7.15, Corollary~10.16 and 
Corollary~10.27 of~\cite{atiyah-macdonald}). It is clear
that $B$ in Definition~\ref{d56} is finitely generated over $R$ if and
only if  $A$ finitely
generated as a $\bfk$-vector space. Therefore, for $A$ finitely
generated over $\bfk$,
free deformations are the same as the flat ones.

The $R$-linearity of deformations implies the following simple
lemma. Recall that all deformations in this sections satisfy~(\ref{Eq87}).

\begin{lemma}
\label{l1}
Let $B = (B,\mu)$ be a deformation as in Definition~\ref{d56}. 
Then the multiplication $\mu$ in $B$ 
is determined by its restriction to $A \otimes A \subset B \otimes
B$. Likewise, every equivalence of deformations $\phi : B' \to B''$ is
determined by its restriction to $A \subset B$.
\end{lemma}

\begin{proof}
By~(\ref{Eq87}), each element of $B$ is a finite sum of
elements of the form $ra$, $r \in R$ and $a \in A$, and
$\mu(ra,sb)=rs\mu(a,b)$ by the $R$-bilinearity of $\mu$ for each $a,b
\in A$ and $r,s \in R$. This proves the first statement. 
The second part of the lemma is equally obvious.
\end{proof}

The following proposition will also be useful.

\begin{proposition}
\label{artindef}
Let $B' = (B',\overline{B}' \stackrel{\alpha'}{\to} A)$ 
and $B'' = (B'',\overline{B}'' \stackrel{\alpha''}{\to} A)$ be 
$R$-deformations of an associative algebra $A$. Assume
that $R$ is either a local Artinian ring or a complete local ring. 
Then every homomorphism $\phi: B'\to B''$  
of $R$-algebras such that $\overline {\phi} =
\alpha''^{-1}\circ \alpha'$ is an equivalence of deformations.
\end{proposition}

\begin{proof}[Sketch of proof]
We must show that $\phi$ is invertible.
One may consider a formal inverse of $\phi$ in the 
form of an expansion in the successive quotients of the maximal
ideal. If $R$ is Artinian, this formal inverse
has in fact only finitely many terms and hence it is an actual
inverse of~$\phi$. If $R$ is complete, this formal expansion is convergent.
\end{proof}

We leave as an exercise to prove that each $R$-deformation of $A$
in the sense of Definition~\ref{d56} is equivalent to a deformation of
the form $(B,\overline{B} \stackrel{\it can}{\to} A)$, with $B = R \ot A$
(equality of $\bfk$-vector spaces) and ${\it can}$ the canonical
map $\overline{B} = \bfk \ot_R( R \ot A) \to A$ given by 
\[
{\it can}(1 \ot_R (1 \ot a)) := a,\ \mbox {for } a \in A.
\]  
Two deformations $(B,\mu')$ and $(B,\mu'')$ of this type are
equivalent if and only if there exists an $R$-algebra isomorphism
$\phi : (B,\mu') \to (B,\mu'')$ which reduces, under the
identification ${\it can} : \overline{B} \to A$, to the
identity $\id_A : A \to A$. Since we will be interested only in equivalence
classes of deformations, we will assume that all deformations are
of the above special form.

\begin{definition}\label{FormalDeformation}
A \emph{formal deformation} is a deformation, in the sense of
Definition~\ref{d56}, over the complete local augmented ring $\bfk[[t]]$. 
\end{definition}

\begin{exercise}
Is $\bfk[x,y,t]/(x^2+txy)$ a  formal deformation of $\bfk[x,y]/(x^2)$?
\end{exercise}

\begin{theorem}
\label{DeformaceRovnice}
A formal deformation $B$ of $A$ is given by a family 
\[
\{\mu_i:A\otimes A\to A\ |\ i\in\mb{N}\}
\] 
satisfying $\mu_0(a,b)=ab$ (the multiplication in $A$) and
\[
\leqno{\mbox {\rm ($D_k$)}}\quad
\sum_{i+j=k,\ i,j\geq 0}\mu_i(\mu_j(a,b),c)=\sum_{i+j=k,\
 i,j\geq 0}\mu_i(a,\mu_j(b,c))\quad\textrm{for all }a,b,c\in A
\]
for each $k \geq 1$.
\end{theorem}

\begin{proof}
By Lemma~\ref{l1}, the multiplication $\mu$ in $B$ is determined by
its restriction to $A \ot A$.
Now expand $\mu(a,b)$, for $a,b \in A$, into the power series
\[
\mu(a,b)=\mu_0(a,b)+t\mu_1(a,b)+t^2\mu_2(a,b)+\cdots
\]
for some $\bfk$-bilinear functions $\mu_i : A \ot A \to A$, $i \geq 0$.
Obviously, $\mu_0$ must be the multiplication in $A$. 
It is easy to see that $\mu$ is
associative if and only if ($D_k$) are satisfied for each $k \geq 1$.
\end{proof}

\begin{remark}\label{HochMu1}
Observe that $(D_1)$ reads 
\[
a\mu_1(b,c)-\mu_1(ab,c)+\mu_1(a,bc)-\mu_1(a,b)c=0
\]
and says precisely that 
$\mu_1\in\Lin(A^{\otimes 2},A)$ is a Hochschild cocycle, 
$\delta_\Hoch(\mu_1)=0$, see Definition~\ref{D5}.
\end{remark}

\begin{example}
\label{Ex56}
Let us denote by $H$ the group
\[
H := \{u = \id_A + \phi_1 t + \phi_2 t^2 + \cdots|\ \phi_i \in \Lin(A,A)\},
\] 
with the multiplication induced by the composition of linear maps. By
Proposition~\ref{artindef}, formal deformations $\mu' =\mu_0 + \mu'_1t
+ \mu'_2t^2 + \cdots$ and $\mu'' =\mu_0 + \mu''_1t + \mu''_2t^2 +
\cdots$ of $\mu_0$ are 
equivalent if and only if  
\begin{equation}
\label{E100}
u \circ (\mu_0 + \mu'_1t + \mu'_2t^2 + \cdots)
= (\mu_0 + \mu''_1t + \mu''_2t^2 + \cdots) \circ (u \ot u).
\end{equation}
  
\end{example}

We close this section by formulating some classical 
statements~\cite{gerstenhaber:AM63,gerstenhaber:Ann.ofMath.64,%
gerstenhaber:Ann.ofMath.66}  
which reveal the connection between deformation theory of associative
algebras and the 
Hochschild cohomology.
As suggested by Remark~\ref{HochMu1}, 
the first natural object to look at is $\mu_1$. This motivates the following

\begin{definition}
\label{idef}
An \emph{infinitesimal deformation} of an algebra $A$ is a 
$D$-deformation of $A$, where
\[
D:=\bfk[t]/(t^2)
\]
is the local Artinian ring of \emph{dual numbers}.
\end{definition}

\begin{remark}
\label{R2}
One can easily prove an analog of Theorem~\ref{DeformaceRovnice} for
infinitesimal deformations, namely that there is a one-to-one
correspondence between infinitesimal deformations of $A$ 
and $\bfk$-linear maps $\mu_1:A\otimes A\to A$
satisfying $(D_1)$, that is, by Remark~\ref{HochMu1}, Hochschild
2-cocycles of $A$ with coefficients in itself. But we can 
formulate a stronger statement:
\end{remark}

\begin{theorem}\label{InfDefHoch}
There is a one-to-one correspondence between the space of equivalence
classes of
infinitesimal deformations of $A$  and the second Hochschild cohomology
$H^2_\Hoch(A,A)$ of $A$ with coefficients in itself.
\end{theorem}

\begin{proof}
Consider two infinitesimal deformations of $A$ given by
multiplications $*'$ and $*''$, respectively.  As we observed in
Remark~\ref{R2}, these deformations are determined by Hochschild 2-cocycles
$\mu_1',\mu_1'' : A \otimes A \to A$, via equations
\begin{align} 
\label{eqstar} 
a*'b&=ab+t\mu'_1(a,b) 
\\
\nonumber 
a*''b&=ab+t\mu_1''(a,b), \quad a,b \in A.
\end{align}
Each equivalence $\phi$ of deformations $*'$ and $*''$  is
determined by a $\bfk$-linear map $\phi_1 : A \to A$,
\begin{gather} 
\label{eqphi}
\phi(a)=a+t\phi_1(a), \quad a\in A,
\end{gather}
the invertibility of such a $\phi$ follows from
Proposition~\ref{artindef} but can easily be checked directly. 
Substituting~\eqref{eqstar} and~\eqref{eqphi} into
\begin{equation}
\label{E5}
\phi(a*'b) = \phi(a)*''\phi(b), \quad a,b \in A,
\end{equation}
one obtains
\[
\phi(ab+t\mu_1'(a,b)) = (a+t\phi_1(a))*''(b+t\phi_1(b))
\]
which can be further expanded into 
\[
ab+t\phi(\mu_1'(a,b)) 
= ab+t(a\phi_1(b))+t(\phi_1(a)b)+t\mu_1''(a+t\phi_1(a),b+t\phi_1(b))
\]
so, finally,
\[
ab+t\mu_1'(a,b) = ab+t(a\phi_1(b)+\phi_1(a)b)+t\mu_1''(a,b).
\]
Comparing the $t$-linear terms, we see that~(\ref{E5}) is equivalent to
\[
\mu_1'(a,b)=\delta_\Hoch\phi_1(a,b)+\mu_1''(a,b).
\]
We conclude 
that infinitesimal deformations given 
by $\mu_1',\mu_1''\in C^2_\Hoch(A,A)$ are equivalent if and only
if they differ by a coboundary, that is, if and only if  
$[\mu_1']=[\mu_1'']$ in $H^2_\Hoch(A,A)$.
\end{proof}

Another classical result is:

\begin{theorem}
\label{rr}
Let $A$ be an associative algebra such that $H_\Hoch^2(A,A)=0$.
Then all formal deformations of $A$ are equivalent to $A$.
\end{theorem}

\begin{proof}[Sketch of proof]
If $*',*''$ are two formal deformations of $A$,
one can, using the assumption $H^2_\Hoch(A,A)=0$, as in the proof of
Theorem~\ref{InfDefHoch} find a $\bfk$-linear map $\phi_1:A\to A$ 
defining an equivalence of $(B,*')$ to $(B,*'')$ 
modulo $t^2$. Repeating this process, 
one ends up with an equivalence $\phi=\id+t \phi_1 +t^2
\phi_2+\cdots$ of formal deformations $*'$ and $*''$.
\end{proof}

\begin{definition}
An \emph{$n$-deformation} of an algebra $A$ is an 
$R$-deformation of $A$ for $R$ the local Artinian ring $\bfk[t]/(t^{n+1})$.
\end{definition}

We have the following version of 
Theorem~\ref{DeformaceRovnice} whose proof is obvious.

\begin{theorem}
An $n$-deformation of $A$ is given by a family 
\[
\{\mu_i:A\otimes A\to A\ |\ 1\leq i\leq n\}
\]
of $\bfk$-linear maps satisfying ($D_k$) of
Theorem~\ref{DeformaceRovnice}  for $1 \leq k \leq n$.
\end{theorem}

\begin{definition}
An $(n+1)$-deformation of $A$ given by $\{\mu_1,\ldots,\mu_{n+1}\}$ 
is called an \emph{extension} of the $n$-deformation given by
$\{\mu_1,\ldots,\mu_n\}$.
\end{definition}

Let us rearrange $(D_{n+1})$ into
\begin{eqnarray*}
\lefteqn{
-a\mu_{n+1}(b,c)+\mu_{n+1}(ab,c)-\mu_{n+1}(a,bc)+\mu_{n+1}(a,b)c=\hskip 2cm}
\\
&& \hskip 2cm
=\sum_{i+j=n+1,\ i,j>0}\left(\mu_i(a,\mu_j(b,c)) -\mu_i(\mu_j(a,b),c)\right)
\end{eqnarray*}
Denote the trilinear function in the right-hand side 
by $\mf{O}_n$ and interpret it as an element of
$C^3_\Hoch(A,A)$,
\begin{equation}
\label{E7}
 \mf{O}_n :=  \sum_{i+j=n+1,\
 i,j>0}\left(\mu_i(a,\mu_j(b,c))- \mu_i(\mu_j(a,b),c)\right) \in
 C^3_\Hoch(A,A). 
\end{equation}
Using the Hochschild differential recalled in Definition~\ref{D5},
one can rewrite ($D_{n+1}$) as
\[
\delta_\Hoch(\mu_{n+1})=\mf{O}_n.
\]
We conclude that,  if an $n$-deformation extends to an 
$(n+1)$-deformation, then $\mf{O}_n$ is a~Hochschild coboundary. 
In fact, one can prove:

\begin{theorem}
{}For any $n$-deformation, the Hochschild cochain 
$\mf{O}_n\in C^3_\Hoch(A,A)$
defined in~(\ref{E7}) is a
cocycle, $\delta_\Hoch(\mf{O}_n)=0$. Moreover, 
$[\mf{O}_n]=0$ in $H^3_\Hoch(A,A)$ if and only if the  
$n$-deformation $\{\mu_1,\ldots,\mu_n\}$ extends into some $(n+1)$-deformation.
\end{theorem}

\begin{proof}
Straightforward. 
\end{proof}

\noindent
{\bf Geometric deformation theory.\/} 
Let us turn our attention back
to the variety  of structure constants $\Asss(V)$ 
recalled in Section~\ref{s2}, page~\pageref{p8}.
Elements of $\Asss(V)$ are associative $\bfk$-linear multiplications 
$\cdot: V\otimes V\to V$ and there is a natural left action $\cdot
\mapsto \cdot_\phi$ of $\GL(V)$ on $\Asss(V)$ given by 
\begin{equation}
\label{E732}
a \cdot_\phi b:=\phi(\phi^{-1}(a)\cdot\phi^{-1}(b)), 
\end{equation}
for $a,b\in V$ and $\phi\in\GL(V)$.
We assume that $V$ is finite dimensional.

\begin{definition}
Let $A$ be an algebra with the underlying vector space $V$ 
interpreted as a point in the variety of structure constants, $A \in Ass(V)$.
The algebra $A$ is called (geometrically) \emph{rigid} if
the $\GL(V)$-orbit of $A$ in $\Asss(V)$ is Zarisky-open.
\end{definition}

Let us remark that, if $\bfk = {\mathbb R}$ or ${\mathbb C}$, then, 
by~\cite[Proposition~17.1]{nijenhuis-richardson:BAMS66}, the
$\GL(V)$-orbit of $A$ in $\Asss(V)$ is Zarisky-open if and only if it
is (classically) open.
The following statement whose proof can be found
in~\cite[\S~5]{nijenhuis-richardson:BAMS66} specifies the relation between the
Hochschild cohomology and geometric rigidity, compare also
Propositions~1 and~2 of~\cite{felix:Bull.Soc.Math.France80}. 

\begin{theorem}
Suppose that the ground field is algebraically closed.
\hfill\break
\hglue\parindent (i) If $H^2_\Hoch(A,A)=0$ then $A$ is rigid, and
\hfill\break
\hglue\parindent
(ii) if $H^3_\Hoch(A,A)=0$ then 
$A$ is rigid if and only if  $H^2_\Hoch(A,A)=0$.
\end{theorem}

\noindent
{\bf Three concepts of rigidity.} One says that an
associative algebra is {\em infinitesimally rigid\/} if 
$A$ has only trivial  (i.e.~equivalent to
$A$) infinitesimal deformations. 
Likewise, $A$ is {\em analytically rigid\/}, if all formal
deformations of $A$ are trivial. 

By Theorem~\ref{InfDefHoch}, $A$ is infinitesimally rigid if and only
if $H^2_\Hoch(A,A) = 0$. Together with Theorem~\ref{rr} this
establishes the first implication in the following display which in
fact holds over fields of arbitrary characteristic
\[
\mbox {infinitesimal rigidity} \Longrightarrow 
\mbox {analytic rigidity} \Longrightarrow 
\mbox {geometric rigidity}.
\]
The second implication in the above display 
is~\cite[Theorem~3.2]{gerstenhaber-schack:JPAA86}. 
Theorem~7.1 of the same paper
then says that in characteristic zero, the analytic and geometric
rigidity are equivalent concepts:
\[
\mbox {analytic rigidity} 
\stackrel{\mbox {\raisebox{.3em}{\scriptsize \rm char.~0}}}\Longleftrightarrow 
\mbox {geometric rigidity}.
\]

\noindent 
{\bf Valued deformations.} The authors of~\cite{goze-remm:JAA04}
studied $R$-deformations of finite-dimensional algebras in the case when
$R$ was a valuation ring~\cite[Chapter~5]{atiyah-macdonald}. 
In particular, they considered deformations over the  non-standard
extension ${\mathbb C}^*$ of the field of complex numbers, and called these
${\mathbb C}^*$-deformations {\em perturbations\/}. They argued,
in~\cite[Theorem~4]{goze-remm:JAA04}, that an algebra $A$ admits only trivial
perturbations if and only if it is geometrically rigid.

\begin{remark}
\label{R8}
An analysis parallel to the one presented in this section 
can be made for any class of ``reasonable''
algebras, where ``reasonable'' are algebras over quadratic Koszul
operads~\cite[Section~II.3.3]{markl-shnider-stasheff:book} 
for which the deformation cohomology is given by
a~``standard construction.'' Let us emphasize that most of
``classical'' types of algebras (Lie, associative, associative
commutative, Poisson, etc.) are ``reasonable.''
See also~\cite{balavoine:CM97,balavoine:jpaa98}.
\end{remark}

\section{Structures of (co)associative (co)algebras}
\label{s4}

Let $V$ be a $\bfk$-vector space.  In this section we recall, in
Theorems~\ref{CorrespAlg} and~\ref{CorrespCoAlg}, the following
important correspondence between (co)algebras and differentials:
\begin{gather}
\nonumber 
\{\textrm{coassociative coalgebra structures on the vector space }V \}
\\
\nonumber 
\label{COAL}
\updownarrow 
\\
\nonumber 
\{\textrm{quadratic differentials on the free 
             associative algebra generated by $V$}\}.
\end{gather}
and its dual version:
\begin{gather}
\nonumber 
\{\textrm{associative algebras on the vector space }V \}
\\
\nonumber 
\label{AL}
\updownarrow 
\\
\nonumber 
\{\textrm{quadratic differentials on the ``cofree'' coassociative 
coalgebra cogenerated   by $V$}\}.
\end{gather}

The reason why we put `cofree' into parentheses will become clear
later in this section.
Similar correspondences exist for any ``reasonable'' (in the sense
explained in Remark~\ref{R8}) class of algebras,
see~\cite[Theorem~8.2]{fox-markl:ContM97}. We will in fact need only
the second correspondence but, since it relies on coderivations of
``cofree'' coalgebras, we decided to start with the first one which is
simpler to explain.

\begin{definition}
\label{d5}
The \emph{free  associative algebra} generated by a vector space $W$ 
is an associative algebra $\mc{A}(W)\in\Ass$ together with a 
linear map $W\to \mc{A}(W)$ having the following property:

For every $A\in\Ass$ and a linear map $W\To^\varphi A$, 
there exists a unique algebra homomorphism $\mc{A}(W)\to A$ 
making the diagram:
\begin{diagram}
W & \rTo & \mc{A}(W) \\
& \rdTo_\varphi & \dDashto \\
& & A \\
\end{diagram}
commutative.
\end{definition}

The free associative
algebra on $W$ is uniquely determined up to isomorphism. An example is
provided by the \emph{tensor algebra} 
$T(W):=\bigoplus_{n=1}^\infty W^{\otimes n}$ with the
inclusion $W = W^{\otimes 1} \hookrightarrow T(W)$.
There is a natural grading on $T(W)$ given by the number of
tensor factors, 
\[
T(W)=\textstyle \bigoplus_{n=0}^\infty T^n(W),
\]
where  $T^n(W) :=
W^{\otimes n}$ for $n \geq 1$ and $T^0(W) := 0$. 
Let us emphasize that the tensor algebra as defined above is {\em
  nonunital\/}, the unital version can be obtained by taking $T^0(W) :
= \bfk$.

\begin{convention}
We are going to consider graded algebraic objects. Our choice of signs 
will be dictated by the principle 
that whenever we commute two
``things'' of degrees $p$ and~$q$, respectively,
we multiply the sign by $(-1)^{pq}$. This rule is sometimes called the
\emph{Koszul sign convention}. As usual, non-graded (classical)
objects will be, when necessary, considered as graded ones
concentrated in degree~$0$.
\end{convention}

Let $f'  : V' \to W'$ and $f'' : V'' \to W''$ be homogeneous maps of
graded vector spaces. The
Koszul sign convention implies that the value of $(f' \ot f'')$ on the
product $v' \ot v'' \in V' \ot V''$ of homogeneous elements equals
\[
(f' \ot f'')(v' \ot v'') := (-1)^{\deg(f'')\deg(v')}f'(v') \ot f''(v'').
\] 
In fact, the Koszul sign convention is determined by the above rule
for evaluation.

\begin{definition}
Assume $V = V^*$ is a graded vector space, $V=\bigoplus_{i\in\mb{Z}}V^i$. The
\emph{suspension operator} $\uparrow$ assigns to $V$ the graded 
vector space $\sssp V$ with $\mb{Z}$-grading $(\sssp
V)^i:=V^{i-1}$. There is a natural degree +1 map 
$\uparrow: V \to \sssp V$ that sends $v \in V$ into its suspended
copy $\sssp v \in \sssp V$.
Likewise, the \emph{desuspension operator} $\downarrow$ changes the
grading of $V$
according to the rule $\hbox{$(\dessp V)^i$}:=V^{i+1}$. 
The corresponding degree $-1$ map $\downarrow : V \to \dessp V$ 
is defined in the obvious~way. The suspension (resp.~the desuspension) of $V$
is sometimes also denoted $s V$ or $V[-1]$ (resp.~$s^{-1} V$ or $V[1]$). 
\end{definition}

\begin{example}
If $V$ is an  un-graded vector space, then $\uparrow V$ 
is $V$ placed in degree $+1$  and $\downarrow V$ is $V$ placed in
degree $-1$. 
\end{example}

\begin{remark}
In the ``superworld'' of $\mb{Z}_2$-graded objects, the operators 
$\uparrow$ and $\downarrow$ agree and 
coincide with the \emph{parity change operator}.
\end{remark}

\begin{exercise}\label{signlemma}
Show that the Koszul sign convention implies
$(\downarrow\otimes\downarrow)\circ(\uparrow\otimes\uparrow)=-\id$ or,
more generally, 
\[
\downarrow^{\otimes n}\circ \uparrow^{\otimes n}= 
\uparrow^{\otimes n}\circ \downarrow^{\otimes n}= (-1)^{\frac
{n(n-1)}2} \id
\]
for an arbitrary $n \geq 1$. 
\end{exercise}

\begin{definition}
A \emph{derivation} of an associative algebra $A$ 
is a linear map $\theta:A\to A$ satisfying the \emph{Leibniz rule} 
\[
\theta(ab)=\theta(a)b+a\theta(b)
\] 
for every $a,b\in A$. Denote $\Der(A)$ the set of all derivations of $A$.
\end{definition}

We will in fact need a graded version of the above definition:

\begin{definition}
\label{d987}
A \emph{degree $d$ derivation} of a $\mb{Z}$-graded algebra 
$A$ is a degree $d$ linear map 
$\theta:A\to A$ satisfying the graded Leibniz rule
\begin{equation}
\label{E8}
\theta(ab)=\theta(a)b+(-1)^{d|a|}a\theta(b)
\end{equation}
for every homogeneous element $a\in A$ of degree $|a|$ and for every $b\in A$. 
We denote $\Der^d(A)$ the set of all degree $d$ derivations of $A$. 
\end{definition}

\begin{exercise}
\label{Ex2}
Let $\mu : A \ot A \to A$ be the multiplication of $A$. 
Prove that~(\ref{E8}) is equivalent~to 
\[
\theta\mu=\mu(\theta\otimes\id)+\mu(\id\otimes\theta).
\]
Observe namely how the signs in the right hand side of~(\ref{E8}) are
dictated by the Koszul convention.
\end{exercise}

\begin{proposition}
\label{T0}
Let $W$ be a graded vector space and $T(W)$ the tensor algebra
generated by $W$ with the induced grading. For any $d$, there is a
natural 
isomorphism
\begin{equation}
\label{EE4}
\Der^d(T(W))\cong \Lin^d(W,T(W)),
\end{equation}
where $\Lin^d(-,-)$ denotes the space of degree $d$ $\bfk$-linear maps.
\end{proposition}

\begin{proof}
Let $\theta\in\Der^d(T(W))$ and $f := \theta|_W : W \to T(W)$. 
The Leibniz rule~(\ref{E8}) implies that, for homogeneous elements
$w_i\in W$, $1 \leq i \leq n$, 
\begin{align*}
\theta(w_1\otimes\ldots\otimes w_n) &
= f(w_1)\otimes w_2\otimes\ldots\otimes w_n 
+ (-1)^{d|w_1|} w_1\otimes f(w_2)\otimes\ldots\otimes w_n + \cdots
\\
&= \sum_{i=1}^n (-1)^{d(|w_1| + \cdots + |w_{i-1}|)}
w_1\otimes\ldots\otimes f(w_i)\otimes\ldots\otimes w_n
\end{align*}
which reveals that $\theta$ is determined by its restriction $f$ on
$W$.  On the other hand, given a degree $d$ linear map $f : W \to
T(W)$, the above formula clearly defines a derivation $\theta \in
\Der^d(T(W))$.  The correspondence
\[
\Der^d(T(W)) \ni \theta \longleftrightarrow f: = 
\theta|_W \in \Lin^d(W,T(W))
\] 
is the required isomorphism~(\ref{EE4}). 
\end{proof}

\begin{exercise}
\label{Ex09}
Let $\theta \in \Der^d(T(W))$, $f := \theta|_V$ and $x \in
T^2(W)$. Prove that
\[
\theta(x) = (f \otimes \id + \id \otimes f)(x).
\]
\end{exercise}

\begin{definition}
A derivation $\theta\in\Der^d(T(W))$ is called \emph{quadratic} 
if $\theta(W)\subset T^2W$. A~degree $1$ derivation $\theta$ is 
a \emph{differential} if
$\theta^2 = 0$. 
\end{definition}

\begin{exercise}
Prove that the isomorphism of Proposition~\ref{T0} restricts to
\[
\Der^d_2(T(W)) \cong \Lin^d(W,T^2(W)),
\]
where $\Der^d_2(T(W))$ is the space of all quadratic degree $d$
derivations of $T(W)$. 
\end{exercise}

\begin{definition}
\label{d77}
Let $V$ be a vector space.
A \emph{coassociative coalgebra} structure on 
$V$ is given by a linear map 
$\Delta:V\to V\otimes V$
satisfying 
\[
(\Delta\otimes \id)\Delta=(\id\otimes\Delta)\Delta
\]
(the \emph{coassociativity}).
\end{definition}

We will need, in Section~\ref{s7}, also a cocommutative version of
coalgebras:

\begin{definition}
A coassociative coalgebra $A = (V,\Delta)$ 
as in Definition~\ref{d77} is {\em cocommutative\/} if 
\[
T \Delta = \Delta
\]
with the swapping map $T : V \ot V \to V \ot V$ 
given by 
\[
T(v' \ot v'') : =
(-1)^{|v'||v''|} v'' \ot v'
\] 
for homogeneous $v',v'' \in V$.  
\end{definition}

\begin{theorem}\label{CorrespAlg}
Let $V$ be a (possibly graded) vector space.  Denote $\Coass(V)$ the
set of all coassociative coalgebra structures on $V$ and
$\Diff^1_2(T(\sssp V))$ the set of all quadratic differentials on
the tensor algebra $T(\sssp V)$.  Then there is a natural isomorphism
\[
\Coass(V)\cong \Diff^1_2( T(\sssp V)).
\]
\end{theorem}

\begin{proof}
Let $\chi\in  \Diff^1_2(T(\sssp V))$.
Put $f:=\chi|_{\uparrow V}$ 
so that $f$ is a degree $+1$ map 
$\sssp V\to \sssp V\otimes\sssp V$.
By Exercise~\ref{Ex09} (with $W := \sssp V$, $\theta := \chi$ 
and $x := f(\sssp v)$), 
\[
0=\chi^2(\sssp v)=\chi(f(\sssp v))=
(f\otimes\id+\id\otimes f)(f(\sssp v))
\]
for every $v \in V$, therefore
\begin{equation}
\label{eq:assoc}
(f\otimes\id+\id\otimes f)f=0.
\end{equation}
We have clearly described a one-to-one correspondence between quadratic
differentials $\chi \in \Diff^1_2(T(\sssp V))$ and degree $+1$
linear maps $f \in \Lin^1(\sssp V, T^2(\sssp V))$
satisfying~(\ref{eq:assoc}).

Given $f: \sssp V \to \sssp V \otimes \sssp V$ as above, 
define the map $\Delta : V \to V \ot V$ by the commutative diagram
\begin{diagram}
\sssp V & \rTo^f & \sssp V\otimes\sssp V \\
\uTo<\susp & & \uTo>{\susp \otimes\susp} \\
V & \rTo^\Delta & V\otimes V\\
\end{diagram}
i.e., by Exercise~\ref{signlemma}, 
\[
\Delta:= (\susp \otimes\susp)^{-1}\circ f\circ\susp 
=- (\desusp\otimes\desusp)\circ f\circ\susp.
\] 
Let us show
that~(\ref{eq:assoc}) is equivalent to the coassociativity of
$\Delta$. We have
\begin{align*}
(\Delta\otimes\id)\Delta 
&= \left(-(\desusp\otimes\desusp)f\susp\otimes\id\right)
\left(-(\desusp\otimes\desusp)f\susp\right)=
\left((\desusp\otimes\desusp)f\susp\otimes\id\right)
(\desusp\otimes\desusp)f\susp
\\
&= ((\desusp\otimes\desusp)f\otimes\desusp)f\susp
=-(\desusp\otimes\desusp\otimes\desusp)(f\otimes\id)f\susp.
\end{align*}
The minus sign in the last term appeared because we interchanged $f$ (a
``thing'' of degree $+1$) with $\desusp$ (a ``thing'' of degree
$-1$). Similarly
\begin{align*}
(\id\otimes\Delta)\Delta &
= \left(\id\otimes(-(\desusp\otimes\desusp))f\susp\right)
\left(-(\desusp\otimes\desusp)f\susp\right)=
 \left(\id\otimes(\desusp\otimes\desusp)f\susp\right)
(\desusp\otimes\desusp)f\susp
\\
&= (\desusp\otimes(\desusp\otimes\desusp)f)f\susp=
(\desusp\otimes\desusp\otimes\desusp)(\id\otimes f)f\susp,
\end{align*}
so~(\ref{eq:assoc}) is indeed equivalent to $(\Delta \ot \id)\Delta =
(\id \ot \Delta)\Delta$. This finishes the proof.
\end{proof}

We are going to dualize Theorem~\ref{CorrespAlg} to get a description of
associative algebras, not {\em co\/}algebras. First, we need a dual
version of the tensor algebra:

\begin{definition}
The underlying vector space $T(W)$ of the tensor algebra
with the comultiplication 
$\Delta:T(W)\to T(W)\otimes T(W)$ defined by 
\[
\Delta(w_1\otimes\ldots\otimes w_n)
:=\sum_{i=1}^{n-1}(w_1\otimes\ldots\otimes w_i)
\otimes(w_{i+1}\otimes\ldots\otimes w_n)
\]
is a coassociative coalgebra denoted $\cT(W)$ 
and called the \emph{tensor coalgebra}.
\end{definition}

\begin{warning}
\label{warn}
Contrary to general belief, the coalgebra $\cT{}(W)$ with the projection
$\cT{}(W) \to W$ \underline{is} \underline{not}
cofree in the category of coassociative
coalgebras! Cofree coalgebras (in the sense of the obvious dual of
Definition~\ref{d5}) are surprisingly complicated
objects~\cite{fox:JPAA93a,smith:TopAppl03,hazewinkel:JPAA03}. 
The coalgebra $\cT{}(W)$ is, however, cofree in the subcategory  of
coaugmented nilpotent 
coalgebras~\cite[Section~II.3.7]{markl-shnider-stasheff:book}. 
This will be enough for our purposes.
\end{warning}

In the following dual version of Definition~\ref{d987} we use
\emph{Sweedler's convention} expressing the comultiplication in a
coalgebra $C$ as $\Delta(c) = \sum c_{(1)} \otimes c_{(2)}$, $c \in
C$.

\begin{definition}
A \emph{degree $d$ coderivation} of a $\mb{Z}$-graded coalgebra
$C$ is a linear degree $d$ map $\theta:C\to C$ satisfying the dual
Leibniz rule
\begin{equation}
\label{E9}
\Delta\theta(c)=\sum\theta(c_{(1)})\otimes c_{(2)}
+\sum(-1)^{d|c_{(1)}|} c_{(1)} \otimes\theta(c_{(2)}),
\end{equation}
for every $c\in C$.  Denote the set of all degree $d$ coderivations of
$C$ by $\CoDer^d(C)$.
\end{definition}

As in Exercise~\ref{Ex2} one easily proves that~(\ref{E9}) is 
equivalent to 
\[
\Delta\theta=(\theta\otimes\id)\Delta+(\id\otimes\theta)\Delta.
\]
Let us prove the dual of Proposition~\ref{T0}:

\begin{proposition}
\label{P23}
Let $W$ be a graded vector space. For any $d$, there is a natural isomorphism
\begin{equation}
\label{EE4bis}
\CoDer^d(\cT{}(W))\cong \Lin^d(T(W),W).
\end{equation}
\end{proposition}

\begin{proof}
{}For $\theta\in\CoDer^d(T(W))$ and $s \geq 1$ denote $f_s \in
\Lin^d(T^s(W),W)$ the composition
\begin{equation}
\label{E354}
f_s :
T^s(W)
\stackrel{\theta|_{T^s(W)}}{-\hskip-.3em-\hskip-.7em\longrightarrow} 
\cT{}(W) \stackrel{{\rm proj.}}{-\hskip-.7em\longrightarrow} W.
\end{equation}
The dual Leibniz rule~(\ref{E9}) implies that, 
for $w_1,\ldots,w_n\in W$ and $n \geq 0$, 
\begin{eqnarray*}
\lefteqn{
\theta(w_1\otimes\ldots\otimes w_n)
:=}
\\
\nonumber 
&&\sum_{s \geq 1}\sum_{i=1}^{n-s+1}
(-1)^{d(|w_1| + \cdots + |w_{i-1}|)}
w_1\otimes\ldots\otimes w_{i-1}\otimes 
f_s(w_i\otimes\ldots\otimes w_{i+s-1})\otimes w_{i+s}
\otimes\ldots\otimes w_n,
\end{eqnarray*}
which shows that $\theta$ is
uniquely determined by $f:= f_0 + f_1 +f_2 + \cdots \in
\Lin^d(T(W),W)$. 
On the other hand, it is easy to verify that for any
map $f \in\Lin^d(T(W),W)$ decomposed into the sum of its
homogeneous components, the above formula defines a coderivation $\theta
\in\CoDer^d(T(W))$. This finishes the proof.
\end{proof}

\begin{definition}
\label{d0986}
The composition $f_s : T^s(W) \to W$ defined in~(\ref{E354}) is called
the $s$th {\em corestriction\/} of the coderivation $\theta$. 
A coderivation $\theta\in\CoDer^d(T(W))$ is \emph{quadratic} if 
its $s$th corestriction is non-zero only for $s=2$. 
A degree~$1$ coderivation $\theta$ is a 
\emph{differential} if $\theta^2 = 0$. 
\end{definition}

Let us finally formulate a dual version of Theorem~\ref{CorrespAlg}.

\begin{theorem}
\label{CorrespCoAlg}
Let $V$ be a graded vector space.
Denote $\CoDiff^1_2(\cT{}(\dessp V))$ the set of all quadratic differentials 
on the tensor coalgebra $\cT{}(\dessp V)$.
One then has a natural isomorphism
\begin{equation}
\label{E66}
\Asss(V) \cong \CoDiff^1_2(\cT{}(\dessp V)).
\end{equation}
\end{theorem}

\begin{proof}
Let $\chi \in \CoDiff^1_2(\cT{}(\dessp V))$ and $f : \dessp V \ot \dessp
V \to \dessp V$ be the 2nd corestriction of $\chi$. Define $\mu : V \ot V \to
V$ by the diagram
\begin{diagram}
\dessp V \ot \dessp V & \rTo^f & \dessp V 
\\
\uTo<{\desusp \ot \desusp} & & \uTo>{\desusp} 
\\
V\ot V & \rTo^\mu & \hphantom{.} V.\\
\end{diagram}
The correspondence $\chi \leftrightarrow \mu$ is then the required
isomorphism. This can be verified by dualizing the steps of the
proof of Theorem~\ref{CorrespAlg} so we can safely leave the details
to the reader.
\end{proof}

\section{dg-Lie algebras and the Maurer-Cartan equation}
\label{s5}

\begin{definition}
\label{d1279}
A graded Lie algebra is a $\mb{Z}$-graded vector space
\[
\mf{g}=\bigoplus_{n\in\mb{Z}}\mf{g}^n
\] 
equipped with a degree $0$
bilinear map $[-,-] : \mf{g} \ot \mf{g} \to \mf{g}$ (the
\emph{bracket}) which is graded antisymmetric,~i.e. 
\begin{equation}
\label{E11}
[a,b]=-(-1)^{|a||b|}[b,a]
\end{equation}
for all homogeneous $a,b \in \mf{g}$, and
satisfies the graded Jacobi identity:
\begin{equation}
\label{E12}
[a,[b,c]]+(-1)^{|a|(|b|+|c|)}[b,[c,a]]+(-1)^{|c|(|a|+|b|)}[c,[a,b]]=0
\end{equation}
for all homogeneous $a,b,c \in \mf{g}$.
\end{definition}

\begin{exercise}
Write the axioms of graded Lie algebras in an element-free form that
would use only the bilinear map $l : = [-,-] : \mf{g} \ot \mf{g}
\to \mf{g}$ and its iterated compositions, and the operator of
``permuting the inputs'' of a multilinear map. Observe how the Koszul
sign convention helps remembering the signs in~(\ref{E11})
and~(\ref{E12}).
\end{exercise}

\begin{definition}
\label{D765}
A \emph{dg-Lie algebra} (an abbreviation for \emph{differential 
graded Lie algebra}) is a graded Lie algebra $L=\bigoplus_{n \in
  {\mathbb Z}}L^n$
as in Definition~\ref{d1279}
together with a degree $1$ linear map $d:L\to L$ which is

-- a degree $1$ derivation , i.e.~$d[a,b]=[da,b]+(-1)^{\abs{a}}[a,db]$
   for homogeneous $a,b \in L$, and

-- a differential, i.e.~$d^2=0$.
\end{definition}

Our next aim is to show that the Hochschild complex 
$\left(C^*_\Hoch(A,A),\delta_\Hoch\right)$ of an
associative algebra recalled in Definition~\ref{D5}
has a natural bracket which turns it into a dg-Lie algebra. We start with
some preparatory material.

\begin{proposition}\label{bracketCoDer}
Let $C$ be a graded coalgebra. For coderivations 
$\theta,\phi\in\CoDer(C)$ define
\[
[\theta,\phi]:=\theta \circ \phi-(-1)^{|\theta||\phi|}\phi\circ \theta.
\] 
The bracket $[-,-]$ makes 
$\CoDer(C)=\bigoplus_{n\in\mb{Z}}\CoDer^n(\mc{C})$ a graded Lie algebra.
\end{proposition}

\begin{proof}
The key observation is that $[\theta,\phi]$ is a coderivation 
(note that {\em neither\/} $\theta\circ \phi$ {\em nor\/} 
$\phi \circ \theta$ are coderivations!).
Verifying this and the properties of a graded Lie bracket is
straightforward and will be omitted.
\end{proof}

\begin{proposition}
\label{differential}
Let $C$ be a graded coalgebra and $\chi\in\CoDer^1(\mc{C})$ such that 
\begin{equation}
\label{E13}
[\chi,\chi]=0,
\end{equation}
where $[-,-]$ is the bracket of Proposition~\ref{bracketCoDer}. Then 
\[
d(\theta):=[\chi,\theta]\qquad \textrm{for }\theta\in\CoDer(\mc{C})
\]
is a differential that makes $\CoDer(\mc{C})$ a dg-Lie algebra.
\end{proposition}

Observe that, since $|\chi|=1$,~(\ref{E13})
does not tautologically follow from the graded antisymmetry~(\ref{E11}).

\begin{proof}[Proof of Proposition~\ref{differential}]
The graded Jacobi identity~(\ref{E12}) implies
that, for each homogeneous~$\theta$, 
\[
[\chi,[\chi,\theta]] 
=-(-1)^{|\theta|+1}[\chi,[\theta,\chi]] - [\theta,[\chi,\chi]].
\]
Now use the graded antisymmetry 
$[\theta,\chi]=(-1)^{|\theta|+1}[\chi,\theta]$ 
and the assumption $[\chi,\chi]=0$ to conclude from the above display that
\[
[\chi,[\chi,\theta]] = -[\chi,[\chi,\theta]],
\]
therefore, since the characteristic of the ground field is zero,
\[
d^2(\theta)=[\chi,[\chi,\theta]]=0,
\]
so $d$ is a differential. The derivation property of $d$ with respect
to the bracket can be verified in the same way and we leave it as an
exercise to the reader.
\end{proof}

In Proposition~\ref{differential} we saw that coderivations of a
graded coalgebra form a dg-Lie algebra. Another example of a dg-Lie
algebra is provided by the Hochschild cochains
of an associative algebra (see Definition~\ref{D5}). We need the following:

\begin{definition}
\label{d90}
For $f\in \Lin(\otexp V{(m+1)},V)$, $g\in \Lin(\otexp V{(n+1)},V)$ and
$1 \leq i \leq m+1$ define $f\circ_i g\in \Lin(\otexp V{(m+n+1)},V)$
by
\[
(f\circ_i g) := 
f\left(\id_V^{\otimes (i-1)} \otimes g \otimes \id_V^{\otimes (m-i+1)}\right).
\]
Define also
\[
f\circ g:=\sum_{i=1}^{m+1}(-1)^{n(i+1)}f\circ_i g
\]
and, finally, 
\[
[f,g]:=f\circ g-(-1)^{mn}g\circ f.
\]
The operation $[-,-]$ is called the \emph{Gerstenhaber bracket} (our
definition however differs from the original one of~\cite{gerstenhaber:AM63}
by the overall sign $(-1)^n$).
\end{definition}

Let $A$ be an associative algebra with the underlying space
$V$. Since, by Definition~\ref{D5}, 
$C^{*+1}_\Hoch(A,A) = \Lin(\otexp V{(*+1)},V)$, the
structure of Definition~\ref{d90} defines a degree $0$ operation
$[-,-] : C^{*+1}_\Hoch(A,A) \otimes C^{*+1}_\Hoch(A,A) \to
C^{*+1}_\Hoch(A,A)$ called again the Gerstenhaber bracket.
We leave as an exercise the proof of

\begin{proposition}
\label{p456}
The Hochschild cochain complex of an associative algebra together with
the Gerstenhaber bracket form a dg-Lie algebra $C^{*+1}_\Hoch(A,A) =
(C^{*+1}_\Hoch(A,A),[-,-],\delta_\Hoch)$.
\end{proposition}

The following theorem gives an alternative description of the dg-Lie algebra of
Proposition~\ref{p456}.  

\begin{theorem}
\label{complexes}
Let $A$ be an associative algebra with multiplication $\mu : V \ot V
\to V$ and $\chi \in \CoDiff^1_2(\cT{}(\dessp V))$ the coderivation
that corresponds to $\mu$ in the correspondence of
Theorem~\ref{CorrespCoAlg}. Let  $d := [\chi,-]$ be 
the differential introduced in Proposition~\ref{differential}.
Then there is a natural isomorphism of dg-Lie algebras
\[
\xi: 
\left( 
C^{(*+1)}_\Hoch (A,A),[-,-],\delta_\Hoch
\right)
\stackrel\cong\longrightarrow \left(\CoDer^*(\cT{}(\dessp V)),[-,-],d\right).
\]
\end{theorem}

\begin{proof}
Given $\phi \in C^{n+1}_\Hoch(A,A)=\Lin(V^{\otimes(n+1)},V)$, let $f :
(\dessp V)^{\otimes (n+1)} \to \dessp V$ be the degree $n$ linear map
defined by the diagram
\begin{diagram}
(\dessp V)^{\otimes(n+1)} & \rTo^f & \dessp V\\
\uTo<{\desusp^{\otimes (n+1)}} & & \uTo>\downarrow \\
V^{\otimes(n+1)} & \rTo^{\phi} &  \hphantom{.}V. 
\end{diagram}
By Proposition~\ref{P23}, there exists a unique coderivation $\theta
\in \CoDer^n(\cT{}(\dessp V))$
whose \hbox{$(n+1)$th} corestriction is $f$ and other corestrictions are
trivial. 

The map $\xi :C^{(*+1)}_\Hoch (A,A) \to
\CoDer^*(\cT{}(\dessp V))$ defined by $\xi(\phi) := \theta$ is clearly
an isomorphism. The verification that $\xi$ commutes with the
differentials and brackets is a straightforward though involved 
exercise on the Koszul sign
convention which we leave to the reader.
\end{proof}

\begin{corollary}
\label{c78}
Let $\mu$ be the multiplication in $A$ interpreted as an element of
$C^2_\Hoch(A,A)$, and $f\in C^*_\Hoch(A,A)$.
Then $\delta_\Hoch(f) = [\mu,f]$.
\end{corollary}

\begin{proof}
The corollary immediately follows from
Theorem~\ref{complexes}. Indeed, because $\xi$ commutes with all the
structures, we have
\[
\delta_\Hoch(f) = \xi^{-1} \xi \delta_\Hoch (f) = \xi^{-1} (d(\xi f)) 
= \xi^{-1} [\chi,\xi f] = [\mu,f].
\]
We however recommend as an exercise to verify the corollary directly,
comparing $[\mu,f]$ to the formula for the Hochschild differential.
\end{proof}

\begin{proposition}
\label{asoc}
A bilinear map $\kappa : V \ot V \to V$ 
defines an associative algebra structure on $V$ if and only if  
$[\kappa,\kappa]=0$.
\end{proposition}

\begin{proof}
By Definition~\ref{d90} (with $m=n=1$),
\[
\frac 12[\kappa,\kappa]=
\frac 12\left(\rule{0pt}{1.2em}
\kappa\circ\kappa-(-1)^{mn}\kappa\circ\kappa\right)
= \kappa\circ\kappa=\kappa\circ_1\kappa-\kappa\circ_2\kappa
= \kappa(\kappa \ot \id_V) - \kappa(\id_V \ot \kappa),
\]
therefore $[\kappa,\kappa]=0$ is indeed equivalent to the
associativity of $\kappa$.
\end{proof}

\begin{proposition}\label{MC}
Let $A$ be an associative algebra with the underlying vector space $V$
and the multiplication $\mu: V \ot V \to V$. Let 
$\nu\in C^2_\Hoch(A,A)$ be a Hochschild 2-cochain.
Then $\mu+\nu \in C^2_\Hoch(A,A) = \Lin(\otexp V2,V)$ 
is associative if and only if  
\begin{equation}
\label{e0}
\delta_\Hoch(\nu)+\frac{1}{2}[\nu,\nu]=0.
\end{equation}
\end{proposition}

\begin{proof}
By Proposition~\ref{asoc}, $\mu+\nu$ is associative if and only if  
\[
0=\frac 12 [\mu+\nu,\mu+\nu]=
\frac 12 \left\{ \rule {0pt}{1.3em}
[\mu,\mu]+[\nu,\nu]+[\mu,\nu]+[\nu,\mu]\right\}  =
\delta_\Hoch (\nu) + \frac 12 [\nu,\nu].
\]
To get the rightmost term, we used the fact that,
since $\mu$ is associative, $[\mu,\mu] = 0$ by Proposition~\ref{asoc}. 
We also observed that  
$[\mu,\nu] = [\nu,\mu]= \delta_\Hoch(\nu)$ by Corollary~\ref{c78}. 
\end{proof}

A bilinear map $\nu : V \ot V \to V$ such that  $\mu+\nu$ is associative
can be viewed as a deformation of $\mu$. This suggests that~(\ref{e0}) is
related to deformations. This is indeed the case, as we will
see later in this section.
Equation~(\ref{e0}) is a particular case of the Maurer-Cartan equation
in a arbitrary dg-Lie algebra:

\begin{definition}\label{defMC1}
Let $L = (L,[-,-],d)$ be a dg-Lie algebra.
A degree~$1$ element 
$s\in L^1$ is \emph{Maurer-Cartan} if it satisfies the 
\emph{Maurer-Cartan equation} 
\begin{equation}
\label{e111}
ds+\frac{1}{2}[s,s]=0.
\end{equation}
\end{definition}

\begin{remark}
The Maurer-Cartan equation (also called the {\em Berikashvili
equation\/}) along with its clones and
generalizations is one of the most important
equations in mathematics. For instance, a version of the Maurer-Cartan
equation describes the differential of a left-invariant form,
see~\cite[I.\S 4]{kobayashi-nomizu}.  
\end{remark}

Let $\mf{g}$ be a dg-Lie algebra over the ground field $\bfk$. 
Consider the dg-Lie algebra $L$ over
the power series ring $\bfk[[t]]$ defined as
\begin{equation}
\label{E143}
L:=\mf{g}\otimes(t),
\end{equation}
where $(t)\subset \bfk[[t]]$ is the ideal generated by $t$.  Degree
$n$ elements
of $L$ are expressions $f_1 t+ f_2 t^2+\cdots$,
$f_i\in\mf{g}^n$ for $i \geq 1$.  
The dg-Lie structure on $L$ is induced from that of $\mf{g}$ in an
obvious manner.
Denote by $\MC(\mf{g})$ the
set of all Maurer-Cartan
elements in $L$.  Clearly, a~degree $1$ element 
$s =f_1 t+ f_2 t^2+\cdots$ is Maurer-Cartan   
if its components $\{f_i \in \mf{g}^1\}_{i \geq 1}$
satisfy the  equation:
\[
\leqno{\mbox {($\it MC_k$)}}
\hskip 5cm df_k + \displaystyle\frac 12 \sum_{i+j = k}[f_i,f_j] = 0 
\]
for each $k \geq 1$.

\begin{example}
\label{Ex6}
Let us apply the above construction to the
Hochschild complex of an associative algebra $A$ with the
multiplication $\mu_0$, that is, take
$\mf{g} := C^{*+1}_\Hoch(A,A)$ with the Gerstenhaber bracket and the
Hochschild differential.  In this case, one easily sees that
$({\it MC}_k)$ for $s =\mu_1t+\mu_2t^2+\cdots$, $\mu_i\in
C^2_\Hoch(A,A)$ is precisely equation $(D_k)$ of
Theorem~\ref{DeformaceRovnice}, $k \geq 1$, compare also calculations on
page~\pageref{E7}. 
We conclude that $\MC(\mf{g})$ is
the set of infinitesimal  deformations of $\mu_0$. 
\end{example}

Let us recall that each Lie algebra $\mf{l}$ can be equipped with a group
structure with the multiplication given by the {\em Hausdorff-Campbell
formula\/}:
\begin{equation}
\label{HC}
x \cdot y := x+y + \frac12[x,y] + \frac1{12}
([x,[x,y]] + [y,[y,x]]) + \cdots 
\end{equation}
assuming a suitable condition that guarantees that the
above infinite sum makes sense in $\mf{l}$, 
see~\cite[I.IV.\S7]{serre:65}. 
We denote $\mf l$ with 
this multiplication by $\exp(\mf{l})$.  Formula~(\ref{HC}) is obtained by
expressing the right hand side of
\[
x \cdot y = \log(\exp(x)\exp(y)),
\]
where 
\[
\exp(a): = 1 + a + \frac 1{2!}a^2 + \frac 1{3!} a^3 + \cdots,\quad
\log(1+a) := a  - \frac12a^2 + \frac13a^3 - \cdots, 
\]
in terms of iterated commutators of non-commutative variables $x$ and $y$.

Using this construction, we introduce the \emph{gauge group} of
$\mf{g}$ as
\[
{\rm G}(\mf{g})
:=\exp(L^0),
\]
where $L^0 = \mf{g}^0 \ot (t)$ is the Lie subalgebra of degree zero
elements in $L$ defined in~(\ref{E143}). Let us fix an element $\chi
\in \mf{g}^1$.  
The gauge group then acts on $L^1
= \mf{g}^1 \ot (t)$ by the formula
\begin{equation}
\label{e788}
x \cdot l := l + [x,\chi + l] + \frac 1{2!} [x,[x,\chi + l]] + \frac
1{3!}[x,[x,[x,\chi + l]]] + \cdots,\ x \in {\rm G}(\mf{g}),  l \in L^1, 
\end{equation}
obtained by expressing the right hand side of
\begin{equation}
\label{E39}
x \cdot l = \exp (x) (\chi + l) \exp (-x) - \chi
\end{equation}
in terms of iterated commutators. Denoting $d \chi := [\chi,\chi]$,
formula~(\ref{e788}) reads
\begin{equation}
\label{E391}
x \cdot l = l + dx + [x,l] + \frac 12 
\left\{\rule{0pt}{1.2em} [x,dx] +  [x,[x,l]]\right\}
+ \frac 13\left\{\rule{0pt}{1.2em} [x,[x,dx]] +  [x,[x,[x,l]]]\right\}
+ \cdots
\end{equation}

\begin{lemma}
\label{l89}
Action~(\ref{E391})
of \ ${\rm G}(\mf{g})$ on $L^1$ preserves the space $\MC(\mf{g})$ of
solutions of the Maurer-Cartan equation.
\end{lemma}

\begin{proof}
We will prove the lemma under the assumption that 
$\mf{g}$ is a dg-Lie algebra whose differential $d$ has the form
$d = [\chi,-]$ for some $\chi \in \mf{g}^1$ satisfying $[\chi,\chi] =
0$ (see Proposition~\ref{differential}). The proof of the general case
is a straightforward, though involved, verification.

It follows from~(\ref{E39}) that 
$\chi + x \cdot l = \exp (x) (\chi + l) \exp (-x)$, 
i.e.~$x$ transforms $\chi + l$ into $\exp
(x) (\chi + l) \exp (-x)$.
Under the assumption $d = [\chi,-]$, the Maurer-Cartan equation for
$l$ is equivalent to $[\chi + l,\chi + l]=0$. The Maurer-Cartan equation
for the transformed $l$ then reads
\[
[\exp (x) (\chi + l) \exp (-x),\exp (x) (\chi + l) \exp (-x)] = 0,
\]
which can be rearranged into
\[
\exp (x) [\chi + l,\chi + l] \exp (-x) = 0.
\] 
This finishes the proof.
\end{proof}

Thanks to Lemma~\ref{l89}, it makes sense to consider
\[
\Def(\mf{g}):=\MC(\mf{g})/\Gg,
\] 
the moduli space of solutions of the Maurer-Cartan equation in $L =
\mf{g} \ot (t)$.

\begin{example}
\label{Ex567}
Let us return to the situation in Example~\ref{Ex6}. In this case
\[
\mf{g}_0= C^1_\Hoch(A,A) = \Lin(A,A),
\] 
with the bracket given by the
commutator of the composition of linear maps.
The~gauge group $\Gg$ consists of elements  
$x =f_1t+f_2t^2+\ldots$, $f_i\in\Lin(A,A)$. It follows from the
definition of the gauge group action that two formal deformations 
$\mu' = 
\mu_0 + \mu'_1t + \mu'_2t^2 + \cdots$ and $\mu'' = 
\mu_0 + \mu''_1t + \mu''_2t^2
+ \cdots$ of $\mu_0$ define the same element in $\Def(\mf{g})$ if
and only if 
\begin{equation}
\label{E99}
\exp(x) (\mu_0 + \mu'_1t + \mu'_2t^2 + \cdots)
= (\mu_0 + \mu''_1t + \mu''_2t^2 + \cdots)(\exp(x) \ot \exp(x))  
\end{equation}
for some $x \in \Gg$. The above formula has an actual,
not only formal, meaning -- all power series make sense because
of the completeness of the ground ring.

On the other hand, recall that in Example~\ref{Ex56}
we introduced the group
\[
H := \{u = \id_A + \phi_1 t + \phi_2 t^2 + \cdots|\ \phi_i \in \Lin(A,A)\}.
\] 
The exponential map $\exp : \Gg \to H$ is a well-defined isomorphism
with the inverse map $\log : H \to \Gg$. We conclude that the equivalence
relation defined by~(\ref{E99}) is the same as the equivalence defined
by~(\ref{E100}) in Example~\ref{Ex56}, therefore $\Def(\mf{g}) =
\MC(\mf{g})/\Gg$ is the moduli space of equivalence classes of formal
deformations of $\mu_0$.
\end{example}

The above analysis can be generalized by replacing, in~(\ref{E143}), $(t)$
by an arbitrary ideal $\mf{m}$ in a local Artinian
ring or in a complete local ring.

\section{$L_\infty$-algebras and the Maurer-Cartan equation}
\label{s7}

We are going to describe a generalization of differential graded Lie
algebras. Let us start by recalling some necessary notions.

Let $W$ be a $\mb{Z}$-graded vector space.
We will denote by
$\ext W$ the \emph{free graded commutative associative 
algebra} over $W$. It is characterized by the
obvious analog of the universal property in Definition~\ref{d5} with
respect to graded commutative associative algebras.
It can be realized as
the tensor algebra $T(W)$ modulo the 
ideal generated by $x\otimes y-(-1)^{\abs{x}\abs{y}}y\otimes
x$. If one decomposes
\[
W=W^\textrm{even}\oplus W^\textrm{odd}
\]
into the even and odd parts, then 
\[
\ext W \cong \bfk[W^\textrm{even}]\otimes E[W^\textrm{odd}],
\]
where the first factor is the polynomial algebra and the second one is
the exterior (Grassmann) algebra.
The algebra $\ext W$ can also be identified with the subspace of
$T(W)$ consisting of graded-symmetric elements (remember we work
over a characteristic zero field).

Denote the product of (homogeneous) elements 
$w_1,\ldots,w_n\in W$ in $\ext W$ by $w_1 \land  \ldots \land
w_n$. For a permutation $\sigma\in\mf{S}_k$ we define the
\emph{Koszul sign} $\varepsilon(\sigma) \in \{-1,+1\}$ by
\[
w_1\wedge\ldots\wedge
w_k=\varepsilon(\sigma)w_{\sigma(1)}\wedge\ldots\wedge
w_{\sigma(k)}
\]
and the \emph{antisymmetric Koszul sign} $\chi(\sigma)\in \{-1,+1\}$ by
\[
\chi(\sigma):=\sgn(\sigma)\varepsilon(\sigma).
\]

\begin{exercise}
Express $\epsilon(\sigma)$ and $\chi(\sigma)$ explicitly in terms of
$\sigma$ and the degrees $|w_1|$,\ldots,$|w_n|$.
\end{exercise}

Finally, a permutation $\sigma\in\mf{S}_n$ is called an 
\emph{$(i,n-i)$-unshuffle} 
if $\sigma(1)<\ldots<\sigma(i)$ and
$\sigma(i+1)<\ldots<\sigma(n)$. The set of all $(i,n-i)$-unshuffles
will be denoted $\mf{S}_{(i,n-i)}$. 

\begin{definition}
\label{Linft}
An $L_\infty$\emph{-algebra} (also called a~\emph{strongly
homotopy Lie} or \emph{sh Lie algebra}) 
is a~graded vector space $V$ together with a system 
\[
l_k:\otimes^kV\to V,\quad k\in\mb{N}
\]
of linear maps of degree $2-k$ subject to the following axioms.

-- Antisymmetry: For every $k\in\mb{N}$, every permutation
   $\sigma\in\mf{S}_k$ and every homogeneous $v_1,\ldots,v_k \in V$,
\begin{equation}
\label{E17}
l_k(v_{\sigma(1)},\ldots,v_{\sigma(k)})=\chi(\sigma)l_k(v_1,\ldots,v_k).
\end{equation}

-- For every $n \geq 1$ and homogeneous  $v_1,\ldots,v_n \in V$,
\[
\leqno{\mbox{\rm ($L_n$)}}\quad  \displaystyle
\sum_{i+j=n+1} (-1)^{i} \sum_{\sigma\in\mf{S}_{i,n-i}}
\chi(\sigma) l_j(l_i(v_{\sigma(1)},\ldots,v_{\sigma(i)}),
v_{\sigma(i+1)},\ldots,v_{\sigma(n)})=0.
\]
\end{definition}

\begin{remark}
\label{kdy_umru}
The sign in~($L_n$) was taken from~\cite{getzler:AM:09}.
With this sign convention, all terms of the (generalized)
Maurer-Cartan equation recalled in~(\ref{e1323}) below have
$+1$-signs. Our sign convention is related to the original one
in~\cite{lada-markl:CommAlg95,lada-stasheff:IJTP93} via the
transformation $l_n \mapsto (-1)^{n+1 \choose 2 } l_n$. We also used
the opposite grading which is better suited for our purposes -- the
operation $l_k$ as introduced
in~\cite{lada-markl:CommAlg95,lada-stasheff:IJTP93} has degree $k-2$.
\end{remark}

Let us expand axioms $(L_n)$ for $n=1,2$ and $3$.

\noindent 
{\bf Case $n=1$.} For $n=1$ $(L_1)$ reduces to 
$l_1(l_1(v))=0$ for every $v\in V$, 
i.e.~$l_1$ is a degree $+1$ differential.

\noindent 
{\bf Case $n=2$.} 
By~(\ref{E17}), $l_2:V\otimes V\to V$ is a linear degree $0$ map which
is graded antisymmetric, 
\[
l_2(v,u)=-(-1)^{\abs{u}\abs{v}}l_2(u,v)
\]
and ($L_n$) for $n=2$ gives
\[
\leqno{\mbox{\rm ($L_2$)}} \hskip 3.5cm
l_1(l_2(u,v))=l_2(l_1(u),v)+(-1)^\abs{u}l_2(u,l_1(v))
\]
meaning that 
$l_1$ is a graded derivation with respect to the multiplication $l_2$.
Writing $d:=l_1$ and $[u,v]:=l_2(u,v)$, 
($L_2$) takes more usual form 
\[
d[u,v]=[du,v]+(-1)^\abs{u}[u,dv].
\]

\noindent 
{\bf Case $n=3$.} 
The degree $-1$ graded antisymmetric map $l_3:\otimes^3V\to V$
satisfies ($L_3$):
\begin{eqnarray*}
\lefteqn{
(-1)^{\abs{u}\abs{w}}[[u,v],w]+
(-1)^{\abs{v}\abs{w}}[[w,u],v]+(-1)^{\abs{u}\abs{v}}[[v,w],u]=}
\\
&&
=(-1)^{\abs{u}\abs{w}}(dl_3(u,v,w)+
l_3(du,v,w)+(-1)^\abs{u}l_3(u,dv,w)+(-1)^{\abs{u}+\abs{v}}l_3(u,v,dw)).
\end{eqnarray*}
One immediately recognizes the three terms of the Jacobi identity in 
the left-hand side and the $d$-boundary of the trilinear map $l_3$ 
in the right-hand side.
We conclude that the bracket $[-,-]$ satisfies the Jacobi 
identity modulo the homotopy $l_3$.

\begin{example}
\label{dgLie}
If all structure operations $l_k$'s
of an $L_\infty$-algebra $L = (V,l_1,l_2,l_3,\ldots)$  
except $l_1$  vanish, then $L$ is just a dg-vector space
with the differential $d=l_1$.  If all $l_k$'s except $l_1$ and $l_2$
vanish, then $L$ is our familiar dg-Lie algebra 
from Definition~\ref{D765} with
$d=l_1$ and the Lie bracket $[-,-]=l_2$. In this sense, dg-Lie 
algebras are particular cases of $L_\infty$-algebras.
\end{example}

\begin{example}
\label{e999}
Let $L' = (V',l'_1,l'_2,l'_3,\ldots)$  and 
$L'' = (V'',l''_1,l''_2,l''_3,\ldots)$ be two
$L_\infty$-algebras. Define their {\em direct sum\/} $L' \oplus L''$ to
be the $L_\infty$-algebra $L' \oplus L''$ with the underlying vector
space $V' \oplus V''$ and structure operations $\{l_k\}_{k \geq 1}$ given by
\[
l_k(v'_1\oplus v''_1,\ldots,v'_k\oplus
v''_k):=l'_k(v'_1,\ldots,v'_k)+l''_k(v''_1,\ldots,v''_k),
\] 
for $v'_1,\ldots,v'_k \in V'$, $v''_1,\ldots,v''_k \in V''$.
\end{example}

For a graded vector space $V$ denote $\bvee_k (V)$ the quotient of
$\bigotimes ^k V$ modulo the subspace spanned by elements
\[
v_1 \ot \cdots \ot v_k 
-  \chi(\sigma)\  v_{\sigma(1)} \ot \cdots \ot v_{\sigma(k)}.
\]
The antisymmetry~(\ref{E17}) implies that the structure operations of
an $L_\infty$ algebra can be interpreted as maps
\[
l_k : \bvee_k (V) \to V,\ k \geq 1.
\] 

We are going to give a description of the set of $L_\infty$-structures
on a given graded vector space in terms of coderivations, 
in the spirit of Theorem~\ref{CorrespCoAlg}. To
this end, we need the following coalgebra which will play the role of $\cT(W)$.

\begin{proposition}
The space $\ext(W)$
with the comultiplication 
$\Delta:\ext(W)\to \ext(W)\otimes \ext(W)$ defined by 
\[
\Delta(w_1\land\ldots\land w_n)
:=\sum_{i=1}^{n-1} \sum_{\sigma\in\mf{S}_{i,n-i}}  \epsilon(\sigma)
(w_{\sigma(1)}\land\ldots\land w_{\sigma(i)})
\otimes(w_{\sigma(i+1)}\land\ldots\land w_{\sigma(n)})
\]
is a graded coassociative cocommutative coalgebra. 
We will denote it $\cext(W)$. 
\end{proposition}

\begin{proof}
A direct verification which we leave to the reader as an exercise.
\end{proof}

For the coalgebra $\cext(W)$, the following analog of
Proposition~\ref{P23} holds.

\begin{proposition}
\label{P224}
Let $W$ be a graded vector space. For any $d$, there is a natural isomorphism
\[
\CoDer^d(\cext(W))\cong \Lin^d(\cext(W),W).
\]
\end{proposition}

We leave the proof to the reader. Observe that the coalgebra
$\cext(W)$ is a direct sum
\[
\cext(W) = \bigoplus_{n \geq 1}\cext^n(W)
\]
of subspaces $\cext^n(W)$ spanned by $w_1 \land \ldots \land
w_n$, for $w_1,\ldots,w_n \in W$. One may define the $s$th {\em
corestriction\/} of a coderivation $\theta \in \CoDer(\cext(W))$ as
the composition
\[
f_s :
\cext^s(W) 
\stackrel{\theta|_{\mbox{$\land$}^s(W)}}{-\hskip-.3em-\hskip-.7em%
\longrightarrow} 
\cext(W) \stackrel{{\rm proj.}}{-\hskip-.7em\longrightarrow} W.
\]
As in Definition~\ref{d0986},
a coderivation $\theta\in\CoDer^d(\cext(W))$ is \emph{quadratic} if 
its $s$th corestriction is non-zero only for $s=2$. 
A \emph{differential} is a degree~$1$ coderivation $\theta$ such that  
$\theta^2 = 0$.

\begin{theorem}
\label{Linfthm}
Denote by $L_\infty(V)$ the set of
all $L_\infty$-algebra structures on a graded vector space $V$ and
$\CoDiff^1(\cext(\dessp V))$ the set of differentials
on $\cext(\dessp V)$.  Then there is a bijection
\[
L_\infty(V)\cong \CoDiff^{1}(\cext(\dessp V)).
\]
\end{theorem}

\begin{proof}
Let $\chi \in \CoDiff^1(\cext(\dessp V))$ and $f_n :\cext^n(\dessp
V) \to \dessp V$ the $n$th corestriction of $\chi$, $n \geq 1$. Define
$\overline{l}_n : \bvee_n (V) \to V$ by the diagram  
\begin{diagram}
\cext_n(\dessp V) & \rTo^{f_n} & \dessp V \\
\uTo<{\bigotimes ^n \hskip -.4em \downarrow} & & \uTo>\downarrow \\
\bvee_n (V) & \rTo^{\overline{l}_n} & \hphantom{.}V.
\end{diagram}
It is then a direct though involved verification that the maps 
\begin{equation}
\label{E91}
l_n :=
(-1)^{n+1 \choose 2}\overline{l}_n
\end{equation} 
define an $L_\infty$-structure
on $V$ and that the correspondence $\chi \leftrightarrow
(l_1,l_2,l_3,\ldots)$ is one-to-one. The reason for the 
sign change in~(\ref{E91}) is explained in~Remark~\ref{kdy_umru}.
\end{proof}

\begin{remark}
By Theorem~\ref{Linfthm}, $L_\infty$-algebras can be alternatively
defined as square-zero differentials on ``cofree'' cocommutative 
coassociative coalgebras (the
reason why we put `cofree' into quotation marks is the same as in
Section~\ref{s4}, see also the warning on page~\pageref{warn}). 
Dual forms of these object, i.e.~square-zero differentials on free
commutative associative algebras, are 
\emph{Sullivan models} that have existed in
rational homotopy theory 
since~1977~\cite{sullivan:Publ.IHES77}. The same objects
appeared as generalizations of Lie algebras 
independently in 1982 in a~remarkable
paper~\cite{dauria-fre:NP82}. As homotopy Lie algebras with a coherent
system of higher homotopies, 
$L_\infty$-algebras were recognized  much
later~\cite{hinich-schechtman:AdvinSovMath1993,lada-stasheff:IJTP93}.
\end{remark}

\begin{exercise}
Show that the isomorphism of Theorem~\ref{Linfthm} 
restricts to the isomorphism
\[
{\it Lie}(V)\cong \CoDiff^{1}_2(\cext(\dessp V))
\]
between the set of Lie algebra structures on $V$ and quadratic
differentials on the coalgebra $\cext(\dessp V)$. This isomorphism
shall be compared to the isomorphism in Theorem~\ref{CorrespCoAlg}.
\end{exercise}

Let us make a digression and see what happens when one allows in
the right hand side of~(\ref{E66}) all, not only quadratic,
differentials. The above material indicates that one should expect a
homotopy version of associative algebras. This is indeed so; one gets
the following objects that appeared in~1963~\cite{stasheff:TAMS63} (but
we use the sign convention of~\cite{markl:JPAA92}).

\begin{definition} 
An ${A_\infty}$\emph{-algebra} (also called a {\em strongly homotopy
associative algebra\/}) is a~graded vector space $V$
together with a system 
\[
\mu_k:V^{\otimes k}\to V,\quad k \geq 1,
\]
of linear maps of degree $k-2$ such that
\[
\leqno{\mbox {\rm ($A_n$)}} \hskip .3em \rule{0pt}{2.5em}
\displaystyle
\sum_{\lambda=0}^{n-1}
\sum_{k=1}^{n-\lambda}(-1)^{k+\lambda+k
\lambda+k(\abs{v_1}+\cdots+\abs{v_\lambda})}
\!\cdot\!\mu_{n-k+1}(v_1,\squeezedldots,v_\lambda,\mu_k(v_{\lambda+1},
\squeezedldots,v_{\lambda+k}),v_{\lambda+k+1},\squeezedldots,v_n)=0
\]
for every $n \geq 1$, $v_1,\ldots,v_n \in V$.
\end{definition}

\begin{figure}[t]
\begin{center}
\setlength{\unitlength}{0.000583in}%
\begin{picture}(3907,3435)(2847,-4288)
\thicklines
\put(3601,-1261){\line( 1,-3){300}}
\put(3901,-2161){\line(-2,-5){200}}
\put(3301,-3661){\line(2,5){300}}
\put(6001,-1261){\line(-1,-3){300}}
\put(5701,-2161){\line( 2,-5){200}}
\put(6301,-3661){\line(-2,5){300}}
\put(4801,-2461){\line(-1,-1){600}}
\put(4201,-3061){\line( 1,-1){600}}
\put(4801,-3661){\line( 1, 1){600}}
\put(5401,-3061){\line(-1, 1){600}}
\put(4801,-2461){\line( 0, 1){1500}}
\put(4801,-961){\line(-4,-1){1200}}
\put(3601,-1261){\line(-1,-2){600}}
\put(3001,-2461){\line( 2,-1){1200}}
\put(5401,-3061){\line( 2, 1){1200}}
\put(6601,-2461){\line(-1, 2){600}}
\put(6001,-1261){\line(-4, 1){1200}}
\put(3001,-2461){\line( 1,-4){300}}
\put(3301,-3661){\line( 5,-2){1500}}
\put(4801,-4261){\line( 0, 1){600}}
\put(3901,-2161){\line( 1, 0){825}}
\put(4876,-2161){\line( 1, 0){825}}
\put(6601,-2461){\line(-1,-4){300}}
\put(6301,-3661){\line(-5,-2){1500}}
\put(4801,-961){\makebox(0,0){$\bullet$}}
\put(6001,-1261){\makebox(0,0){$\bullet$}}
\put(6601,-2461){\makebox(0,0){$\bullet$}}
\put(6301,-3661){\makebox(0,0){$\bullet$}}
\put(4801,-4261){\makebox(0,0){$\bullet$}}
\put(4801,-3661){\makebox(0,0){$\bullet$}}
\put(5401,-3061){\makebox(0,0){$\bullet$}}
\put(4801,-2449){\makebox(0,0){$\bullet$}}
\put(4201,-3061){\makebox(0,0){$\bullet$}}
\put(3601,-1261){\makebox(0,0){$\bullet$}}
\put(3001,-2461){\makebox(0,0){$\bullet$}}
\put(3301,-3661){\makebox(0,0){$\bullet$}}
\put(3901,-2161){\makebox(0,0){$\bullet$}}
\put(5701,-2161){\makebox(0,0){$\bullet$}}
\end{picture}
\end{center}
\caption{Saito's portrait of $K_5$.\label{Masahito}}
\end{figure}
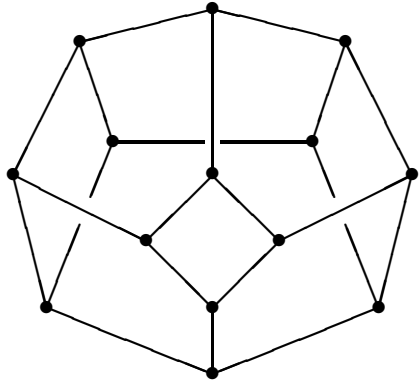

One easily sees
that ($A_1$) means  that $\partial:=\mu_1$ is a degree $-1$ differential,
($A_2$) that the bilinear product $\mu_2 : V \ot V \to V$ commutes
with $\partial$ and ($A_3$) that $\mu_2$ is associative up to the homotopy
$\mu_3$. $A_\infty$-algebras can also be described as algebras over
the cellular chain complex of the
non-$\Sigma$ operad $K = \{K_n\}_{n \geq 1}$ whose $n$th piece is the
$(n-2)$-dimensional convex polytope $K_n$ called the {\em Stasheff
  associahedron\/}~\cite[Section~II.1.6]{markl-shnider-stasheff:book}.   
Let us mention at least that $K_2$ is the point, 
$K_3$ the closed interval and  $K_4$ is the pentagon from Mac Lane's
theory of monoidal categories~\cite{maclane:RiceUniv.Studies63}.
A portrait of $K_5$ due to Masahico Saito is in Figure~\ref{Masahito}.

\begin{theorem}
For a  graded vector space $V$ denote $A_\infty(V)$ the set of all
$A_\infty$-algebra structures on $V$ and $\CoDiff^{1}(\cT{}(\dessp
V))$ the set of all differentials on $\cT{}(\dessp V)$.  Then there is
a natural bijection
\[
A_\infty(V) \cong  \CoDiff^{1}(\cT{}(\dessp V)).
\]
\end{theorem}

\begin{proof}
The isomorphism in the above theorem is of the same nature as the
isomorphism of Theorem~\ref{Linfthm}, but it also involves the `flip'
of degrees since we defined, following~\cite{markl:JPAA92}, 
$A_\infty$-algebras in such a way that the
differential $\partial = \mu_1$ has degree $-1$. We leave the details to
the reader.  
\end{proof}

Let us return to the main theme of this section. Our next task will be to
introduce morphisms of $L_\infty$-algebras.
We start with a simple-minded definition.

Suppose $L' = (V',l'_1,l'_2,l'_3,\ldots)$ and $L'' = 
(V'',l''_1,l''_2,l''_3,\ldots)$ are two $L_\infty$-algebras. A
\emph{strict morphism} is a degree zero linear map $f: V'\to
V''$ which commutes with all structure operations, that is
\[
f(l'_k(v_1,\ldots,v_k))=l''_k(f(v_1),\ldots,f(v_k)),
\]
for each $v_1,\ldots,v_k \in V'$, $k \geq 1$. 

For our purposes we need, however, a subtler notion of morphisms.
We give a definition that involves the isomorphism of Theorem~\ref{Linfthm}. 

\begin{definition}
\label{D456}
Let $L'$ and $L''$ be $L_\infty$-algebras represented 
by dg-coalgebras $(\cext(\dessp
V'),\delta')$ and $(\cext(\dessp
V''),\delta'')$. A {\em (weak) morphism\/} of $L_\infty$-algebras is then 
a morphism of dg-coalgebras $F : (\cext(\dessp
V'),\delta')   \to (\cext(\dessp V''),\delta'')$. 
\end{definition}

Definition~\ref{D456} can be unwrapped. Let $F_k : \cext_k(\dessp
V') \to \dessp V''$ be, for each $k \geq 1$, the composition
\[
\cext^k(\dessp V') 
\stackrel{F}{-\hskip-.3em-\hskip-.7em%
\longrightarrow} 
\cext(\dessp V'') 
\stackrel{{\rm proj.}}{-\hskip-.7em\longrightarrow}\dessp V''.
\]
Define the maps $f_k : \bvee_k V' \to V''$ by the diagram
\begin{diagram}
\cext_k(\dessp V') & \rTo^{F_k} & \dessp V'' \\
\uTo<{\bigotimes^k \hskip -.3em \downarrow} & & \uTo>\downarrow \\
\bvee_k V' & \rTo^{f_k} & \hphantom{.''}V.''
\end{diagram}
Clearly, $f_k$ is a degree $1-k$ linear map. The fact that $F$ is a
dg-morphism can be expressed via a sequence of
axioms ($M_n$), $n \geq 1$, where
($M_n$) postulates the vanishing of 
a combination of $n$-multilinear maps on $V'$ with values in $V''$ 
involving $f_i,l'_i$ and $l''_i$ for
$i \leq n$.

We are not going to write ($M_n$)'s here. 
Explicit axioms for $L_\infty$-maps can be found in~\cite{kajiura-stasheff}, 
see also~\cite[Definition~5.2]{lada-markl:CommAlg95} where
the particular case when $L''$ is a dg-Lie algebra
($l''_k = 0$ for $k \geq 3$) is discussed in detail.
The reader is however encouraged to verify that ($M_1$) says that  
$f_1 : (V',l'_1) \to (V'',l''_1)$ 
is a chain map and that ($M_2$) means that $f_1$ 
commutes with the brackets $l'_2$ and $l''_2$ modulo the
homotopy $f_2$. 

Morphisms of $L_\infty$-algebras $L'$ and $L''$ with underlying vector
spaces $V'$ and $V''$ can
therefore be equivalently defined as systems $f = \{f_k : \bigotimes^k V' \to
V''\}_{k \geq 1}$, where $f_k$ is a degree $1-k$  graded antisymmetric 
linear map, and axioms ($M_n$), $n \geq 1$, are satisfied.
Let us denote by $\catLinfty$ the category of
$L_\infty$-algebras and their morphisms in the sense of
Definition~\ref{D456}.

\begin{exercise}
Show that the category ${\sf strL}_\infty$ of $L_\infty$-algebras and their
strict morphisms can be identified with the (non-full) 
subcategory of ${\sf L}_\infty$ with the same objects and morphisms 
$f = (f_1,f_2,\ldots)$ such that $f_k = 0$ for $k \geq 2$.

Show that the obvious imbedding ${\sf dgLie} 
\hookrightarrow {\sf  L}_\infty$ is not full. This means that there are
more morphisms between dg-Lie algebras considered as
elements of the category ${\sf  L}_\infty$ than in the category of
${\sf dgLie}$. Observe finally that the forgetful
functor 
$\square:\Linf\to\ms{dgVect}$ given by forgetting all structure
operations is not faithful.
\end{exercise}

\section{Homotopy invariance of the Maurer-Cartan equation}
\label{s8}

Let us start with recalling some necessary definitions.

\begin{definition}
A morphism
$f = (f_1,f_2,\ldots) : L' =(V',l'_1,l'_2,\ldots)  
\to L'' =(V'',l''_1,l''_2,\ldots)$ of $L_\infty$-algebras 
is a \emph{weak equivalence} if the chain map
$f_1:(V',l'_1) \to (V'',l_1'')$ induces an isomorphism of cohomology.
\end{definition}

\begin{definition}
An $L_\infty$-algebra $L = (V,l_1,l_2,\ldots)$ is {\em minimal\/} if
$l_1= 0$. It is {\em contractible\/} if $l_k = 0$ for $k \geq
2$ and if $H^*(V,l_1)= 0$. 
\end{definition}

\begin{proposition}
\label{lemma}
A weak equivalence of minimal $L_\infty$-algebras is an isomorphism.
\end{proposition}

\begin{proof}
Let $f = (f_1,f_2,\ldots) : L' \to L''$ be a weak equivalence of
$L_\infty$-algebras.
It follows from the minimality of $L'$ and
$L''$  that the linear part $f_1$ is an isomorphism, thus
the corresponding map $F : (\cext(\dessp
V'),\delta')   \to (\cext(\dessp V''),\delta'')$ induces an
isomorphism of cogenerators. 
It can be easily shown that 
such maps can be inverted.
\end{proof}

The following theorem, which can be found in~\cite{kontsevich:LMPh03},
uses the direct sum of $L_\infty$-algebras recalled in Example~\ref{e999}.

\begin{theorem}
\label{T122}
Each $L_\infty$-algebra is the direct sum of a minimal
and a contractible $L_\infty$-algebra.  
\end{theorem}

Let $L \cong L_m \oplus L_c$ be a decomposition of an
$L_\infty$-algebra $L$
into a minimal $L_\infty$-algebra $L_m$ and a contractible
$L_\infty$-algebra $L_c$. Since the
inclusion $\iota : L_m \to  L_m \oplus L_c \cong L$ is a~weak
equivalence, Theorem~\ref{T122} implies:

\begin{corollary}
\label{c99}
Each $L_\infty$-algebra is weakly equivalent to a minimal one.
\end{corollary}

Corollary~\ref{c99} can also be derived from homotopy invariance
properties of strongly homotopy algebras proved in~\cite{markl:ha}.
Suppose we are 
given an $L_\infty$-algebra $L = (V,l_1,l_2,\ldots)$.
In characteristic zero, two cochain
complexes have the same cochain homotopy type if and only if they have
isomorphic cohomology. In particular, the cochain
complex $(V,l_1)$ is homotopy equivalent to the cohomology
$H^*(V,l_1)$ considered as a complex with trivial differential. 
Move {\bf(M1)} on page~133 of~\cite{markl:ha} now implies that 
there exists an induced minimal $L_\infty$-structure on  
$H^*(V,l_1)$, weakly equivalent to  $L$. Let us remark that an
$A_\infty$-version of Corollary~\ref{c99} was known to Kadeishvili
already in 1985, see~\cite{kadeishvili:ttmi85}.

Remarkably, each $L_\infty$-algebra is, under some mild
assumptions, weakly equivalent to a~dg-Lie algebra. This can be proved
as follows. Suppose $L$ is an  $L_\infty$-algebra represented by a
dg-coalgebra $(\cext (\dessp V),\delta)$. The bar construction
$B(\cext (\dessp V),\delta)$ is a dg-Lie algebra and one may show, 
under an assumption that guarantees the convergence of a spectral
sequence, that $B(\cext (\dessp V),\delta)$ is weakly equivalent to $L$ in
the category of $L_\infty$-algebras. This property is an algebraic
analog of the {\em rectification principle\/} for
$W{\mathcal P}$-spaces provided by the $M$-construction of Boardman
and Vogt, see~\cite[Theorem~II.2.9]{markl-shnider-stasheff:book}.

Let $\mf{g}$ be an $L_\infty$-algebra over the ground field $\bfk$,
with the underlying $\bfk$-vector space $V$. 
Then $V \otimes(t)$, 
where $(t)\subset \bfk[[t]]$ is the ideal generated by $t$, has a
natural induced $L_\infty$-structure. Denote this $L_\infty$-algebra
by $L := \mf{g} \ot (t) = (V \otimes(t), l_1,l_2,l_3,\ldots)$. 
Let $\MC(\mf{g})$ be the
set of all degree $+1$ elements $s \in L^1$ satisfying the {\em generalized
Maurer-Cartan equation\/}
\begin{equation}
\label{e1323}
l_1(s)+ \frac12 l_2(s,s)+ \frac1{3!} l_3(s,s,s) +\cdots+ 
\frac 1{n!}l_n(s,\ldots,s)+\cdots=0.
\end{equation}
When $\mf{g}$ is a dg-Lie algebra, one recognizes the
ordinary Maurer-Cartan equation~(\ref{e111}).

At this moment one  needs to introduce a suitable {\em gauge equivalence\/}
between solutions of~(\ref{e1323}) generalizing the action of
the gauge group $\Gg$ recalled in~(\ref{e788}). Since in
applications of Section~\ref{s9} 
all relevant $L_\infty$-algebras
are in fact dg-Lie algebras, we are not going to describe this
generalized gauge equivalence here, and only refer
to~\cite{kontsevich:LMPh03} instead. We denote
$\Def(\mf{g})$ the set of gauge equivalence classes of solutions
of~(\ref{e1323}). Let us, however, mention that there are examples, as
bialgebras treated in~\cite{markl:ib},
where deformations are described by a fully-fledged
$L_\infty$-algebra.

\begin{example}
\label{1point}
For $\mf{g}$ contractible, $\mf{Def}(\mf{g})$ is the one-point set
consisting of the equivalence class of the trivial solution of~(\ref{e1323}).
Indeed, \[
\MC(\mf{g})=\{s=s_1t+s_2t^2+\ldots\ |\ ds_1 = ds_2 = \cdots = 0\}
\] 
so, by acyclicity, $s_i=db_i$ for some $b_i\in\mf{g}^0$, $i \geq 1$. 
Formula~(\ref{E391}) (with $x = -b_1t_1 - b_2 t_2 - \cdots$ and $l =
s_1t + s_2 t^2 + \cdots$) gives
\[
(-b_1t_1 - b_2 t_2 - \cdots)\cdot (s_1t + s_2 t^2 + \cdots) = 0,
\]
therefore $s = s_1t + s_2 t^2 + \cdots$ is equivalent to the trivial solution.
\end{example}

\begin{example}
Let $\mf{g}'$ and  $\mf{g}''$ be two $L_\infty$-algebras. Then, for
the direct product,
\[
\Def(\mf{g}' \oplus \mf{g}'') 
\cong\Def(\mf{g}')\times\Def(\mf{g}'').
\]
Indeed, it follows from definition that $\MC(\mf{g}' \oplus \mf{g}'') 
\cong\MC(\mf{g}')\times\MC(\mf{g}'')$. This factorization is
preserved by the gauge equivalence. 
\end{example}

The central statement of this section reads:

\begin{theorem}
\label{main}
The assignment $\mf{g} \mapsto \Def(\mf{g})$ extends to a covariant
functor from the category of $L_\infty$-algebras and
their weak morphisms to the category of sets. 
A weak equivalence $f : \mf{g}' \to \mf{g}''$
induces an isomorphism $\Def(f) : \Def(\mf{g}') \cong
\Def(\mf{g}'')$.
\end{theorem}

The above theorem implies that the deformation functor $\Def$ descends
to the localization ${\sf hoL}_\infty$ obtained by inverting weak
equivalences in ${\sf L}_\infty$. By Quillen's
theory~\cite{quillen:LNM43}, ${\sf hoL}_\infty$ is equivalent to the
category of minimal $L_\infty$-algebras and homotopy classes (in an
appropriate sense) of their maps.
This explains the meaning of homotopy
invariance in the title of this section.

\begin{proof}[Proof of Theorem~\ref{main}]
For an $L_\infty$-morphism 
$f = (f_1,f_2,f_3,\ldots) : \mf{g}' \to \mf{g}''$ define
$\MC(f):\MC(\mf{g}')\to \MC(\mf{g}'')$ by
\[
\MC(f)(s):=f_1(s)+\frac 12 f_2(s,s)+\cdots+\frac 1{n!}f_n(s,\ldots,s)+\cdots
\] 
It can be shown that $\MC(f)$ is a well-defined map that descends to
the quotients by the gauge equivalence, giving rise to a map
$\Def(f):\Def(\mf{g}')\to \Def(\mf{g}'')$. 

Assume that $f :\mf{g}' \to \mf{g}''$ above is a weak equivalence.  
By Theorem~\ref{T122}, $\mf{g}'$ decomposes as
$\mf{g}'=\mf{g}'_m\oplus\mf{g}'_c$, with
$\mf{g}'_m$ minimal and $\mf{g}'_c$ contractible, and there is a
similar decomposition $\mf{g}''=\mf{g}''_m\oplus\mf{g}''_c$ for
$\mf{g}''$.  
Define the map $\overline{f} : \mf{g}'_m \to \mf{g}''_m$ by the
commutativity of the diagram
\begin{diagram}
\mf{g}'_m\oplus\mf{g}'_c & \lInto^{i} & \mf{g}'_m \\
\dTo<f & & \dTo>{\overline{f}} \\
\mf{g}''_m\oplus\mf{g}''_c & \rOnto^{p} &\mf{g}''_m
\end{diagram}
in which $i$ is the natural inclusion and $p$ the natural projection.
Observe that $\overline{f}$ is a weak equivalence so it is,
by Proposition~\ref{lemma}, an isomorphism. Therefore, in the
following induced diagram, 
the map $\Def(\overline{f})$ is an isomorphism, too:
\begin{diagram}
\mf{Def}(\mf{g}'_m)\times\mf{Def}(\mf{g}'_c) & \lInto^{\mf{Def}(i)} &
\mf{Def}(\mf{g}'_m) \\ \dTo<{\mf{Def}(f)} & & \dTo>{\mf{Def}(\overline{f})} \\
\mf{Def}(\mf{g}''_m)\times\mf{Def}(\mf{g}''_c) & \rOnto^{\mf{Def}(p)} &
\hphantom{.}\mf{Def}(\mf{g}''_m).
\end{diagram}
Since, by Example \ref{1point}, both $\mf{Def}(\mf{g}'_c)$ and
$\mf{Def}(\mf{g}''_c)$ are points, the maps $\mf{Def}(i)$ and
$\mf{Def}(p)$ are isomorphisms.  We finish the proof by 
concluding that $\Def(f)$ is also an
isomorphism.
\end{proof}

\section{Deformation quantization of Poisson manifolds}
\label{s9}

In this section we indicate the main ideas of  Kontsevich's proof 
of the existence of a~deformation quantization of Poisson
manifolds. Our exposition follows~\cite{kontsevich:LMPh03}. Let us
recall some necessary notions.

\begin{definition}
A {\em Poisson algebra\/} is a vector space
$V$ with operations $\ \cdot : V \ot V \to V$ 
and $\{-,-\}: V \ot V \to V$ such that:

-- $(V,\ \cdot\ )$ is an associative commutative algebra,

-- $(V,\{-,-\})$ is a Lie algebra, and 

-- the map $v\mapsto\{u,v\}$ is a $\cdot$\ -derivation 
for any $u\in V$, i.e. $\{u,v\cdot w\}=\{u,v\}\cdot w+v\cdot  \{u,w\}$.
\end{definition}

\begin{exercise}
Show that
Poisson algebras can be equivalently defined 
as structures with only one operation $\minibullet: V
 \ot V \to V$ such that 
\[ 
u \minibullet (v \minibullet w)=(u \minibullet v)\minibullet w -\frac{1}{3} \!
\left\{(u \minibullet w)\minibullet v+(v \minibullet w)\minibullet u
-(v \minibullet 
u)\minibullet w-(w \minibullet u)\minibullet v)\rule{0pt}{1.2em}\right   \}, 
\]
for each $u,v,w \in V$, see~\cite[Example~2]{markl-remm}.
\end{exercise}

Poisson algebras are `classical limits' of associative 
deformations of commutative associative 
algebras. By this we mean the following. 
Let $A = (V,\ \cdot\ )$ be an associative algebra 
with multiplication $a,b \mapsto
a \cdot b$. 
Consider a formal deformation $(\bfk[[t]]\ot V,\star)$ of $A$ given, as in
Theorem~\ref{DeformaceRovnice},  by a family $\{\mu_i :
A \ot A \to A\}_{i \geq 1}$ by the formula
\begin{equation}
\label{e98}
a\star b:=a\cdot b+t\mu_1(a,b)+t^2\mu_2(a,b)+ t^3 \mu_3(a,b) +\cdots
\end{equation}
for $a,b \in V$. We have the following:

\begin{proposition}
\label{p76}
Suppose $A = (V,\ \cdot\ )$ is a commutative associative
algebra. Then, for an~associative deformation~(\ref{e98}) of $A$, 
\[
\{a,b\}:=\mu_1(a,b)-\mu_1(b,a),\ a,b \in V,
\]
is a Lie bracket such that $P_\star :=  (V,\cdot,\{-,-\})$ 
is Poisson algebra. 
\end{proposition}

\begin{definition}
In the above situation, $P_\star$ is called 
the {\em classical limit\/} of the $\star$-product and 
$(\bfk[[t]] \ot V,\star)$ 
a {\em deformation quantization\/} of the Poisson algebra $P_\star$.
\end{definition}

\begin{proof}[Proof of Proposition~\ref{p76}]
Let us prove first that $\{-,-\}$ is a Lie bracket. The antisymmetry
of $\{-,-\}$ is obvious, 
one thus only needs to verify the Jacobi identity. It is
a standard fact that the  antisymmetrization of an associative multiplication
is a Lie product~\cite[Chapter~I]{serre:65}, 
therefore $[-,-]$ defined by $[x,y] := x \star y - y\star x$
for $x,y \in \bfk[[t]]\otimes A$, is a Lie bracket on $\bfk[[t]]\otimes A$.
We conclude by observing that the Jacobi identity for $\{-,-\}$
evaluated at $a,b,c\in A$ is the term at $t^2$ of the Jacobi identify
for $[-,-]$ evaluated at the same elements.

It remains to verify the derivation property. It is clearly
equivalent to 
\begin{equation}
\label{smrt_se_blizi}
\mu_1(ab,c)-\mu_1(c,ab)-a\mu_1(b,c)+a\mu_1(c,b)-\mu_1(a,c)b+\mu_1(c,a)b
= 0
\end{equation}
where we, for brevity, omitted the symbol for the $\cdot$\ -product.
In Remark~\ref{HochMu1} we observed that $\mu_1$ is a Hochschild
cocycle, therefore 
\[
\rho(a,b,c) :=
a\mu_1(b,c)-\mu_1(ab,c)+\mu_1(a,bc)-\mu_1(a,b)c=0.
\]
A straightforward verification involving the commutativity 
of the $\cdot$\ -product shows that the left hand
side of~(\ref{smrt_se_blizi}) equals $-\rho(a,b,c) + \rho(a,c,b) -
\rho(c,a,b)$. This finishes the proof. 
\end{proof}

Let us recall geometric versions of the above notions.

\begin{definition}
A \emph{Poisson manifold} is a smooth 
manifold $M$ equipped with a Lie product $\{-,-\}: C^\infty(M)
\ot C^\infty(M) \to C^\infty(M)$ on the space of smooth functions 
such that $(C^\infty(M),\ \cdot\ ,\{-,-\})$, where $\cdot$ is 
the standard pointwise multiplication, is a Poisson algebra.
\end{definition}

Poisson manifolds generalize symplectic ones in that the bracket
$\{-,-\}$ need not be induced by a nondegenerate $2$-form.
The following notion was introduced and physically justified 
in~\cite{bayen78}.

\begin{definition}
A \emph{deformation quantization} (also called a {\em star
  product\/}) of a Poisson manifold $M$ is a deformation 
quantization of the Poisson algebra $(C^\infty(M),\cdot,\{-,-\})$ such
that all $\mu_i$'s in~(\ref{e98}) are differential operators.
\end{definition}

\begin{theorem}[Kontsevich~\cite{kontsevich:LMPh03}]
Every Poisson manifold admits a deformation quantization.
\end{theorem}

\begin{proof}[Sketch of Proof]
Maxim 
Kontsevich proved this theorem in two steps. He proved first a~`local' version
assuming $M = {\mathbb R}^{d}$, and then he 
globalized the result to an arbitrary $M$ using ideas of formal
geometry and the language of superconnections.
We are going to sketch only the first step of Kontsevich's proof. 

The idea was to construct two weakly equivalent 
$L_\infty$-algebras $\mf{g}'$, $\mf{g}''$ such that
$\mf{Def}(\mf{g}')$ contained the moduli space of Poisson structures
on $M$ and $\mf{Def}(\mf{g}'')$ was the moduli space of star
products, and then apply Theorem~\ref{main}. In fact, $\mf{g}'$ will
turn out to be an ordinary graded Lie algebra and $\mf{g}''$ a dg-Lie algebra.

\noindent 
{\it-- Construction of $\mf{g}'$.} 
It is the graded Lie algebra of {\em polyvector fields\/} with
the Shouten-Nijenhuis bracket. In more detail,
$\mf{g}'=\bigoplus_{n\geq 0} \mf{g}'^n$ with
\[
\mf{g}'^n :=\Gamma(M,\ext^{n+1}TM), \ n \geq 1,
\]
where $\Gamma(M,\ext^{n+1}TM)$ denotes the space of smooth
sections of the $(n+1)$th exterior power of the tangent bundle $TM$.
The bracket is determined  by
\begin{eqnarray*}
\lefteqn{
[\xi_0\wedge\ldots\wedge\xi_k,\eta_0\wedge\ldots\wedge\eta_l]
:=}
\\
&&:=
\sum_{i=0}^k\sum_{j=0}^l(-1)^{i+j+k}[\xi_i,\eta_j]
\wedge\xi_0\wedge\ldots \wedge\hat{\xi}_i
\wedge\ldots\wedge\xi_k\wedge
\eta_0\wedge\ldots\wedge\hat{\eta}_j\wedge\ldots\wedge\eta_l,
\end{eqnarray*}
where $\xi_1,\ldots,\xi_k,\eta_1,\ldots,\eta_l \in \Gamma(M,TM)$ are
vector fields, $\hat{}$~indicates the omission  
and $[\xi_i,\eta_j]$ in the right hand side denotes the
classical Lie bracket of vector fields $\xi_i$ and $\eta_j$~\cite[I.\S
1]{kobayashi-nomizu}.

Recall that Poisson structures on $M$ are in one-to-one
correspondence with smooth sections $\alpha\in\Gamma(M,\ext^2TM)$ 
satisfying $[\alpha,\alpha]=0$. The corresponding bracket 
of smooth functions $f,g \in C^\infty(M)$ is given by
$\{f,g\} = \alpha(f \otimes g)$.
Since $\mf{g}'$ is just a graded Lie algebra, 
\[
\MC(\mf{g}')=\{s=s_1t+s_2t^2+\ldots \in \mf{g}'^1 \ot (t)\ 
|\ [s,s]=0\}
\] 
therefore clearly $s := \alpha t \in \MC(\mf{g}')$ for each $\alpha
\in \Gamma(M,\ext^2TM)$ defining a Poisson structure. We see that
$\Def(\mf{g}')$ contains the
moduli space of Poisson structures on $M$. 

\vskip .5em

\noindent 
{\it-- Construction of $\mf{g}''$.} It is the dg Lie algebra of
polydiffenential operators, 
\[
\mf{g}''=\bigoplus_{n\geq 0}D^n_{\textrm{poly}}(M),
\]
where 
\[
D^n_{\textrm{poly}}(M)\subset
C^{n+1}_{\Hoch}(C^{\infty}(M),C^{\infty}(M))
\] 
consists of
Hochschild cochains (Definition~\ref{D5}) of the algebra
$C^{\infty}(M)$ given by polydifferential operators. It is clear 
that $D^*_{\textrm{poly}}(M)$ is closed under 
the Hochschild differential and the
Gerstenhaber bracket, so the dg-Lie structure of
Proposition~\ref{p456} restricts to a dg-Lie structure on $\mf{g}''$.
The analysis of Example~\ref{Ex567} shows that
$\Def(\mf{g}'')$ represents equivalence classes of star products. 

\vskip .5em

\noindent 
{\it-- The weak equivalence.} Consider the map 
$f_1:\mf{g'}\to\mf{g''}$ defined by
\[
f_1(\xi_0,\ldots,\xi_k)(g_0,\ldots,g_k)
:=\frac1{(k+1)!} \sum_{\sigma\in\mf{S}_{k+1}}
\sgn(\sigma)\prod_{i=0}^k\xi_{\sigma(i)}(g_i),
\]
for $\xi_0,\ldots,\xi_k\in\Gamma(M,TM)$ and 
$g_0,\ldots,g_k\in C^\infty(M)$.
It is easy to show that $f_1: \hbox{$(\mf{g}',d=0)$} \to
(\mf{g}'',\delta_\Hoch)$ is a chain map. Moreover, a version of
the Kostant-Hochschild-Rosenberg theorem for smooth manifolds proved
in~\cite{kontsevich:LMPh03} states that 
$f_1$ is a cohomology isomorphism.
Unfortunately, $f_1$ {\em does not\/} commute with
brackets. The following central statement
of Kontsevich's approach to deformation quantization says that $f_1$
is, however, the linear part of an $L_\infty$-map:

\begin{formality-theorem}
The map $f_1$ extends to an $L_\infty$-homomorphism $f = (f_1,f_2,f_3,\ldots):
\mf{g}' \to \mf{g}''$.
\end{formality-theorem}

The formality theorem  implies that  
$\mf{g}'$ and $\mf{g}''$ are weakly equivalent in the
category of $L_\infty$-algebras. In other words, the dg-Lie
algebra of polydifferential operators is weakly equivalent to its
cohomology. The `formality' in the name of the theorem is
justified by rational homotopy theory where formal algebras are
algebras having the homotopy type of their cohomology.

Kontsevich's construction of higher $f_i$'s involves coefficients
given as integrals over
compactifications of certain configuration spaces. 
An independent approach of Tamarkin~\cite{tamarkin:deligne} 
based entirely on homological
algebra uses a 
solution of the Deligne conjecture, see also an overview~\cite{hinich:FM03} containing references
to original sources.

\end{proof}

\def\cprime{$'$}

\end{document}